\documentclass[12pt]{amsart}
\usepackage{amsmath,amsthm,enumerate,graphics,amsfonts,latexsym,amsopn,verbatim,amscd,amssymb}

\usepackage{amsmath,amssymb,graphicx,amsthm}
\usepackage{pgfplots}
\usepackage{caption}
\usepackage{geometry}
\usepackage{tikz}

\usepackage{color}
\usepackage{tikz}
\usepackage{amsmath,amssymb,graphicx,tikz,amsthm}
\usetikzlibrary{arrows,matrix}

\theoremstyle{plain}
\newtheorem{thm}{Theorem}[section]

\newtheorem{cor}[thm]{Corollary}

\newtheorem{problem}{Theorem}[section] 
\newtheorem{prob}[problem]{Problem}

\theoremstyle{definition}

\DeclareMathOperator{\cnlspec}{CNL-spec}
\DeclareMathOperator{\cnqspec}{CNSL-spec}
\DeclareMathOperator{\CNRS}{CNRS}
\DeclareMathOperator{\CNL}{CNL}
\DeclareMathOperator{\CNSL}{CNSL}
\DeclareMathOperator{\CN}{CN}
\DeclareMathOperator{\cnrs}{CNRS}
\DeclareMathOperator{\tr}{tr}
\DeclareMathOperator{\Ecn}{E_{CN}}

\tikzset{every path/.style=thick,
	acteur/.style={
		circle,
		fill=red,
		thick,
		inner sep=1pt,
		minimum size=0.15cm
}}

\begin{document}

	\title[Common neighborhood (signless) Laplacian   spectrum and energy]{Common neighborhood (signless) Laplacian  spectrum and energy of CCC-graph}

	\author[F. E. Jannat and R. K. Nath ]{Firdous Ee Jannat and Rajat  Kanti Nath$^*$ }
	\address{Firdous Ee Jannat, Department of Mathematical Science, Tezpur  University, Napaam -784028, Sonitpur, Assam, India.}
	
	\email{firdusej@gmail.com}
	\address{Rajat  Kanti Nath, Department of Mathematical Science, Tezpur  University, Napaam -784028, Sonitpur, Assam, India.}
	
	\email{rajatkantinath@yahoo.com}
\thanks{Corresponding author}		
\begin{abstract} 
In this paper, we consider commuting conjugacy class graph (abbreviated as CCC-graph) of a finite group $G$ which is a graph with vertex set $\{x^G : x \in G \setminus Z(G)\}$ (where $x^G$ denotes the conjugacy class containing $x$) and two distinct vertices $x^G$ and $y^G$ are joined by an edge if there exist some elements $x'\in x^G$ and $y'\in y^G$ such that they commute. We compute common neighborhood (signless) Laplacian spectrum and energy of CCC-graph of finite non-abelian groups whose   central quotient is isomorphic to either $\mathbb{Z}_p \times \mathbb{Z}_p$ (where $p$ is any prime) or the dihedral group $D_{2n}$ ($n \geq 3$); and determine whether CCC-graphs of these groups are common neighborhood (signless) Laplacian hyperenergetic/borderenergetic.	As a consequence, we characterize certain finite non-abelian groups viz. $D_{2n}$, $T_{4n}$, $U_{6n}$, $U_{(n, m)}$, $SD_{8n}$ and $V_{8n}$ such that their CCC-graphs are common neighborhood (signless) Laplacian hyperenergetic/borderenergetic.

\end{abstract}
	
\thanks{ }
	
\maketitle
		
\section{Introduction} \label{S:intro}

Characterizing  finite groups through various graphs defined on them have been an active area of research over the last 50 years. A number of graphs have been defined on groups \cite{PC-21}. Among those, in our paper, we   consider commuting conjugacy class graph (abbreviated as CCC-graph) of a finite  non-abelian group  $G$. For any element $x \in G$, we write $x^G$ to denote the conjugacy class of $G$ containing $x$. The CCC-graph of $G$,  denoted by $\Gamma_G$, is defined as a simple undirected graph whose vertex set is the set of conjugacy classes of non-central elements of $G$ and two vertices $x^G$ and $y^G$  are adjacent if there exists some elements $x'\in x^G$ and $y'\in y^G$ such that $x'y'=y'x'$.
In 2009, Herzog et al. \cite{HLM} introduced the concept of CCC-graph of a group. In 2016, Mohammadian et al. \cite{MEFW-2016} have characterized finite groups such that their CCC-graph is triangle-free. Later on Salahshour and Ashrafi \cite{SA-2020,SA-CA-2020}, obtained structures of CCC-graph of several families of finite CA-groups.  Salahshour \cite{Salah-2020} also described $\Gamma_G$ for the groups whose central quotient is isomorphic to a dihedral group. Characterizations of various classes of finite non-abelian groups through energy, (signless) Laplacian  energy, common neighborhood energy (abbreviated as CN-energy) and genus of their CCC-graphs can be found in \cite{BN-2021,BN-2022,JN-2025}.

The energy of a graph $\mathcal{G}$ (denoted by $E(\mathcal{G})$) is the sum of absolute values of all the eigenvalues of the adjacency matrix of $\mathcal{G}$. This notion was also used in obtaining $\pi$-electron energy of a conjugated carbon molecule in theoretical chemistry. 
The study of energy of a graph was initiated by Gutman \cite{Gutman78} in 1978.     After a long time, in 2006 and 2008, Gutman et al. introduced two more graph-energy-like quantity, known as Laplacian energy \cite{zhou} (denoted by $LE(\mathcal{G})$) and signless Laplacian energy \cite{E-LE-Gutman} (denoted by $LE^+(\mathcal{G})$),  using Laplacian and signless Laplacian eigenvalues of a graph respectively.  
These energies are used to study various properties of graphs (see  \cite{DM-2014, GP-2017}). Later on,  mathematicians have introduced several kinds of graph energies (see \cite{Gutman-Furtula-2017,GUTMAN-2019*}) and studied graph properties. In 2011, Alwardi et al. \cite{ASG}  have introduced CN-energy of a graph. 
Let $\mathcal{G}$ be a graph with vertex set $V(\mathcal{G}) = \{ v_1, v_2, v_3, \cdots, v_n \}$.
Let $C(v_i,v_j)$ be the set of vertices of a graph $\mathcal{G}$ other than $v_i$ and $v_j$ which are adjacent to both $v_i$ and $v_j$. Then the common neighborhood matrix (CN-matrix) of $\mathcal{G}$, denoted by $CN(\mathcal{G})$,  is a matrix of size $n$  whose $(i, j)$-th entry is given by
\[
\CN(\mathcal{G})_{i, j} = \begin{cases} 
	|C(v_i,v_j)|, & \text{ if } i\neq j\\
	0, & \text{ otherwise. } 
\end{cases}
\]
The CN-energy of  $\mathcal{G}$ is the sum of absolute values of all the  eigenvalues of $CN(\mathcal{G})$.
Motivated by the study of  (signless) Laplacian energy, Jannat et al. \cite{FR-2021} have introduced the notions of common neighborhood Laplacian energy (CNL-energy) and common neighborhood signless Laplacian energy (CNSL-energy) of a graph.

The common neighborhood Laplacian matrix (CNL-matrix) and the common neighborhood signless Laplacian matrix (CNSL-matrix) of  $\mathcal{G}$, denoted by $\CNL(\mathcal{G})$ and $\CNSL(\mathcal{G})$, respectively, are given by
\[
\CNL(\mathcal{G}):= \CNRS(\mathcal{G}) - \CN(\mathcal{G}) \text{ and } \CNSL(\mathcal{G}):= \CNRS(\mathcal{G}) + \CN(\mathcal{G}),
\] 
where $\CNRS(\mathcal{G})$ is  a matrix of size $|V(\mathcal{G})| =n$ 
whose $(i, j)$-th entry is  given by
\[
\CNRS(\mathcal{G})_{i, j} =  \begin{cases}
	\underset{k = 1}{\overset{n}{\sum}}\CN(\mathcal{G})_{i, k}, &   \text{ if }  i = j  \text{ and } i = 1, 2, \dots, n\\
	0, & \text{ if } i \ne j.  
\end{cases}
\]

The common neighborhood Laplacian spectrum  of $\mathcal{G}$ (abbreviated as CNL-spectrum and denoted by $\cnlspec(\mathcal{G})$) is the set of eigenvalues of $\CNL(\mathcal{G})$ with multiplicities. We write $\cnlspec(\mathcal{G}) = \{(\alpha_1)^{a_1}, \,(\alpha_2)^{a_2}, \dots,\,(\alpha_k)^{a_k}\}$, where $\alpha_1,\,\alpha_2, \dots, \,\alpha_k$ are the distinct eigenvalues of $\CNL(\mathcal{G})$ with corresponding multiplicities  $a_1,\,a_2, \dots,\,a_k$. Similarly,  common neighborhood signless Laplacian spectrum  of $\mathcal{G}$ (abbreviated as CNSL-spectrum and denoted by $\cnqspec(\mathcal{G})$) is the set of eigenvalues of $\CNSL(\mathcal{G})$ with multiplicities. We write $\cnqspec(\mathcal{G})= \{(\beta_1)^{b_1},\,(\beta_2)^{b_2}, \dots,\,(\beta_{\ell})^{b_{\ell}}\}$, where $\beta_1,\,\beta_2, \dots,\,\beta_{\ell}$ are the distinct eigenvalues of $\CNSL(\mathcal{G})$ with corresponding multiplicities  $b_1,\,b_2, \dots,\,b_{\ell}$. A graph $\mathcal{G}$ is called CNL-integral (CNSL-integral)  if CNL-spectrum  (CNSL-spectrum) contains only integers. The notions of CNL-integral and CNSL-integral graphs were introduced in  \cite{FR-2021} motivated by the notions of integral (introduced by Harary and Schwenk \cite{B10}), L-integral (introduced by Grone and Merris \cite{Grone-94}), Q-integral (introduced by Simic and Stanic \cite{Simic-08}) and CN-integral (introduced by Alwardi et al. \cite{ASG}) graphs. A finite graph is called super integral if it is integral,  L-integral and Q-integral (see \cite{BN-2021}).
Integral graphs have some interests for designing the network topology of perfect state transfer networks (see \cite{E1} and the references there in).

The  CNL-energy  and  CNSL-energy of $\mathcal{G}$, denoted by $LE_{CN}(\mathcal{G})$ and $LE^+_{CN}(\mathcal{G})$ respectively, are defined as 
\begin{equation}\label{LEcn}
	LE_{CN}(\mathcal{G}) := \sum_{i=1}^{k} a_i\left|\alpha_i - \Delta(\mathcal{G}) \right|
\end{equation}
and
\begin{equation}\label{LE+cn}
	LE^+_{CN}(\mathcal{G}) := \sum_{i=1}^{\ell} b_i\left| \beta_i - \Delta(\mathcal{G}) \right|,
\end{equation}
where $\Delta(\mathcal{G}) = \frac{\tr(\CNRS(\mathcal{G}))}{|V(\mathcal{G})|}$ and $\tr(\CNRS(\mathcal{G}))$ is the trace of $\CNRS(\mathcal{G})$. In \cite{FR-2021}, various facets of the CNL-spectrum, CNL-energy, CNSL-spectrum and CNSL-energy of graphs were discussed;  their connections with other well-known graph energies and Zagreb indices were also established. 
It was  observed  that 
\begin{equation}\label{LEcn-Kn}
	LE_{CN}(K_n) = 	LE^+_{CN}(K_n) = 2(n -1)(n - 2),
\end{equation}
where $K_n$ is the complete graph of order $n$. A graph $\mathcal{G}$ of order $n$ is called CNL-hyperenergetic or    CNSL-hyperenergetic according as  $LE_{CN}(\mathcal{G}) > 2(n -1)(n - 2)$ or  $LE^+_{CN}(\mathcal{G}) > 2(n -1)(n - 2)$.   
Further, it is called CNL-borderenergetic (CNSL-borderenergetic) if CNL-energy (CNSL-energy) of $\mathcal{G}$ is equal to  $2(n -1)(n - 2)$. These classes of graphs were introduced in \cite{JN-24,FR-2021} motivated by the notions of various types of hyperenergetic graphs (see \cite{ASG,  FSN-2020,  GLXGF-2015,Gutman1999, TH-2018, Tura-2017, Walikar}).

In Section 2, we compute  CNL-spectrum,  CNSL-spectrum,  CNL-energy  and CNSL-energy of CCC-graphs of finite non-abelian groups whose   central quotient is isomorphic to either $\mathbb{Z}_p \times \mathbb{Z}_p$ (where $p$ is any prime) or the dihedral group $D_{2n}$ ($n \geq 3$). In Section 3, we determine whether CCC-graphs of these groups are CNL-integral, CNSL-integral, CNL-hyperenergetic,   CNSL-hyperenergetic, CNL-borderenergetic and CNSL-borderenergetic. As a consequence, we characterize  the groups  viz. $D_{2n}$, $T_{4n}$, $U_{6n}$, $U_{(n, m)}$, $SD_{8n}$ and $V_{8n}$ such that their CCC-graphs have above mentioned properties. In Subsection 3.1, we compare various CN-energies of  CCC-graphs of the groups $G$ considered in Section 2 and show that $\Ecn(\Gamma_G)  = LE_{CN}(\Gamma_G) = LE^+_{CN}(\Gamma_G)$ or $\Ecn(\Gamma_G) < LE^+_{CN}(\Gamma_G) < LE_{CN}(\Gamma_G)$. We also characterize  the groups  viz. $D_{2n}$, $T_{4n}$, $U_{6n}$, $U_{(n, m)}$, $SD_{8n}$ and $V_{8n}$ such that their CCC-graphs satisfy above mentioned equality/inequality. Finally, we conclude the paper in Section 4 by listing certain problems that arise naturally after our investigation.

\section{Computations of spectrum and energies}
In this section, we compute CNL-spectrum, CNSL-spectrum and their respective energies   of CCC-graphs of various families of non-abelian finite groups. In particular, we consider finite non-abelian groups whose central quotients are isomorphic to $\mathbb{Z}_p \times \mathbb{Z}_p$ (where $p$ is any prime) or $D_{2n}=\langle x,y:~x^{n}=y^2=1,~yxy^{-1}=x^{-1} \rangle$ (for $n \geq 3$) in the following subsections. We shall also consider a generalization of dihedral groups, namely $U_{(n,m)}=\langle x,y:~x^{2n}=y^m=1,~x^{-1}yx=y^{-1}\rangle$ (for  $n\geq 2$, $m\geq 2$) along with the dicyclic group $T_{4n}=\langle x,y:~x^{2n}=1,x^n=y^2,y^{-1}xy=x^{-1} \rangle$ (for $n\geq 2$), the semidihedral groups $SD_{8n}=\langle x,y:~x^{4n}=y^2=1,yxy=x^{2n-1}\rangle$ (for $n\geq 2$),  the groups $U_{6n}=\langle x,y:~x^{2n}=y^3=1,~x^{-1}yx=y^{-1}\rangle$ (for $n\geq 2$)  and 
  $V_{8n}=\langle x,y:~x^{2n}=y^4=1,~yx=x^{-1}y^{-1},~y^{-1}x=x^{-1}y \rangle$ (for $n\geq 2$). The following result is useful in our computation.
\begin{thm}\cite{FR-2021}\label{CN-LE-CNSL-LE^+_Kn}
	Let $\mathcal{G} = l_1 K_{m_1} \cup l_2 K_{m_2}\cup l_3 K_{m_3}$, where $l_iK_{m_i}$ denotes the disjoint union of $l_i$ copies of $K_{m_i}$  for $i = 1, 2, 3$. Then
\begin{align*}
&\cnlspec(\mathcal{G})\\
&\qquad =\left\lbrace (0)^{l_1 + l_2 + l_3}, (m_1(m_1 - 2))^{l_1(m_1 - 1)}, (m_2(m_2 - 2))^{l_2(m_2 - 1)},
(m_3(m_3 - 2))^{l_3(m_3 - 1)}\right\rbrace
\end{align*}  
and
\begin{align*}
\cnqspec(\mathcal{G})=&\left\lbrace (2(m_1 - 1)(m_1 - 2))^{l_1}, ((m_1 - 2)^2)^{l_1(m_1 - 1)}, (2(m_2 - 1)(m_2 - 2))^{l_2},\right.\\
&\left. ((m_2 - 2)^2)^{l_2(m_2 - 1)}, (2(m_3 - 1)(m_3 - 2))^{l_3}, ((m_3 - 2)^2)^{l_3(m_3 - 1)} \right\rbrace.
\end{align*}
\end{thm}

\subsection{Groups whose central quotient is isomorphic to $\mathbb{Z}_p \times \mathbb{Z}_p$}

This class of groups have been considered by Salahshour and Ashrafi \cite[Theorem 3.1]{SA-2020} and showed that 
 \begin{equation}\label{eq-ccc-Zp}
 \Gamma_G = (p + 1)K_{\frac{(p - 1)|Z(G)|}{p}}.
\end{equation} 
 CN-energy of CCC-graphs of this class of groups have been studied in \cite{JN-2025}.  It is worth mentioning that commuting and non-commuting graphs of this class of groups are also studied in \cite{ DBN-2020,PJN-2018,DN-2017, DN-2018, PDN-2018, B17, FN-2020,  NFDS-2021}. In the following theorem, we derive CNL-spectrum, CNSL-spectrum, CNL-energy  and CNSL-energy of CCC-graphs of this class of groups.   
\begin{thm}\label{Z_pZ_p theorem}
	Let $G$ be a finite non-abelian group with  $|Z(G)|=z\geq 2$ and $\frac{G}{Z(G)} \cong \mathbb{Z}_p \times \mathbb{Z}_p$, where $p$ is a prime. Then CNL-spectrum, CNSL-spectrum, CNL-energy and CNSL-energy of CCC-graph of $G$ are given by
	
	$\cnlspec(\Gamma_G)= \left\{(0)^{p+1}, \left(\frac{1}{p^2}(pz - z)(pz - z - 2p)\right)^{\frac{(p + 1)}{p}(pz - z - p)}\right\}$,
	\begin{align*}
		&\cnqspec(\Gamma_G) =\\
		&\qquad\qquad \left\{ \left(\frac{2}{p^2}(pz - z - p)(pz - z - 2p)\right)^{p + 1}, \left(\frac{1}{p^2}(pz - z - 2p)^2\right)^{\frac{(p + 1)}{p}(pz - z - p)} \right\}
	\end{align*}
	and
	\begin{align*}
		LE_{CN}(\Gamma_G) = LE^+_{CN}(\Gamma_G) = \begin{cases}
			\frac{3}{2}, &\text{ for } p=2 \,\&\, z=3\vspace{0.2 cm}\\
			\frac{4 (p-2) (p+1)}{p^2}, &\text{ for } p \geq 2 \,\&\, z=2\\
			\frac{2 (p+1) (p (z-2)-z) (p (z-1)-z)}{p^2}, &\text{ otherwise.}
		\end{cases} 
	\end{align*}
\end{thm}
\begin{proof}
	From \eqref{eq-ccc-Zp}, we have
	$
	\Gamma_G = (p + 1)K_n,
	$
	where $n = \frac{(p - 1)z}{p}$. Therefore, by Theorem \ref{CN-LE-CNSL-LE^+_Kn}, we get
	
	\centerline{
	$\cnlspec(\Gamma_G)= \{(0)^{p+1}, (\frac{1}{p^2}(pz - z)(pz - z - 2p))^{\frac{(p + 1)}{p}(pz - z - p)}\}$ and}
	$\cnqspec(\Gamma_G)= \{ (\frac{2}{p^2}(pz - z - p)(pz - z - 2p))^{p + 1}, (\frac{1}{p^2}(pz - z - 2p)^2)^{(p + 1)\frac{1}{p}(pz - z - p)} \}$.
	\par Here $|V(\Gamma_G)| = \frac{(p^2 - 1)z}{p}$ and $\tr(\cnrs(\Gamma_G)) = \frac{(p - 1) (p + 1) z (p (z - 2) - z) (p (z - 1) - z)}{p^3}$. Therefore, $\Delta(\Gamma_G) = \frac{(p (z - 2) - z) (p (z - 1) - z)}{p^2}$.
	 
Now we calculate CNL-energy of $\Gamma_G$. 	We have
	\begin{align*}
		L_1 &:= \left| 0 - \Delta(\Gamma_G) \right| = \left| - \frac{(p (z - 2) - z) (p (z - 1) - z)}{p^2} \right|.
	\end{align*}
	Let $\alpha_1(p, z) = - (p (z - 2) - z) (p (z - 1) - z)$. Then  $\alpha_1(p, z)= - 2p^2 - 3pz - z^2 + \frac{1}{2}p^2 z(6 - z) + \frac{1}{2}pz^2(4 - p) < 0$ for $p \geq 4$ and $z \geq 6$. It can be seen that $\alpha_1(2, z) = - z(z - 6) - 8 = 1$ or $ \leq 0$ according as $z=3$ or $z\ne3$; $\alpha_1(3, z) = - 2 (z - 3) (2 z - 3) = 2$  or $ \leq 0$ according as $z=2$ or $z\ne2$; $\alpha_1(p, 2) = 2p - 4 \geq 0$; $\alpha_1(p, 3) = - 2p^2 + 9p - 9 = 1$ or $ \leq 0$ according as $p=2$ or $p\ne2$;
	$\alpha_1(p, 4) = - 6 p^2 + 20 p - 16 \leq 0 $ and $\alpha_1(p, 5) = - 12 p^2 + 35 p - 25 \leq 0$. Therefore
	\[
	L_1 = \begin{cases}
		- \frac{(p (z - 2) - z) (p (z - 1) - z)}{p^2}, &\text{ for } p = 2 \,\&\, z=  3; p \geq 2 \, \& \, z =2 \\
		\frac{(p (z - 2) - z) (p (z - 1) - z)}{p^2}, &\text{ otherwise. }
	\end{cases}
	\]
Also
	\begin{align*}
		L_2 &:= \left| \frac{(p z - z) (p z - 2 p - z)}{p^2} - \Delta(\Gamma_G) \right| = \left| \frac{p z - 2 p - z}{p} \right|.
	\end{align*}
	Let $\alpha_2(p,z) = p z - 2 p - z$. Then $\alpha_2(p,z) = (z - 2)p - z  \geq z - 4 \geq 0$ for all $z \geq 4$   since $p \geq 2$.
	It can be seen that  $\alpha_2(p, 2) = -2 < 0$ and $\alpha_2(p, 3) = p - 3 \geq 0$ or $< 0$ according as $p\geq 3$ or $ p = 2$.
	Therefore
	\[ 
	L_2 = \begin{cases}
		- \frac{p z - 2 p - z}{p}, &\text{ for } p=2 \,\&\, z=3; p \geq 2 \,\&\, z=2\\
		\frac{p z - 2 p - z}{p}, &\text{ otherwise.}
	\end{cases}
	\]
	Hence, by \eqref{LEcn}, we get
	\begin{align*}
		LE_{CN}(\Gamma_G) &= (p + 1)\times L_1 + \frac{p + 1}{p}(pz - z - p)\times L_2\\
		&= \begin{cases}
			\frac{3}{2}, &\text{ for }p=2 \,\&\, z=3\vspace{0.2 cm}\\
			\frac{4 (p-2) (p+1)}{p^2}, &\text{ for } p \geq 2 \,\&\, z=2\\ 
			\frac{2 (p+1) (p (z-2)-z) (p (z-1)-z)}{p^2}, &\text{ otherwise. }
		\end{cases}  
	\end{align*}

For CNSL-energy of $\Gamma_G$ we have 
	\begin{align*}
		B_1 :=& \left| \frac{2 (p z - p - z) (p z - 2 p - z)}{p^2}  - \Delta(\Gamma_G) \right| =\left| \frac{(p z-2 p-z) (p z-p-z)}{p^2} \right|\\
		=& - L_1= \begin{cases}
			- \frac{(p z-2 p-z) (p z-p-z)}{p^2}, &\text{ for }p = 2 \,\&\, z=3; p \geq 2 \,\&\, z=2\vspace{0.2 cm}\\
			\frac{(p z-2 p-z) (p z-p-z)}{p^2}, &\text{ otherwise. }
		\end{cases}
	\end{align*}
	Also
	\begin{align*}
		B_2 := \left| \frac{(p z-2 p-z)^2}{p^2} - \Delta(\Gamma_G) \right| &=\left| \frac{- pz + 2 p + z}{p} \right|= - L_2\\
		&= \begin{cases}
			\frac{2p + z - pz}{p}, &\text{ for } p=2 \,\&\, z=3; p \geq 2 \,\&\, z=2\vspace{0.1 cm}\\
			- \frac{2p + z - pz}{p}, &\text{ otherwise. }
		\end{cases}
	\end{align*}
	Hence, by \eqref{LE+cn}, we get
	\begin{align*}
		LE^+_{CN}(\Gamma_G) &= (p + 1)\times B_1 + \frac{p + 1}{p}(pz - z - p)\times B_2\\
		&= \begin{cases}
			\frac{3}{2}, &\text{ for } p=2 \,\&\, z=3\vspace{0.1 cm}\\
			\frac{4 (p-2) (p+1)}{p^2}, &\text{ for } p \geq 2 \,\&\, z=2 \\
			\frac{2 (p+1) (p (z-2)-z) (p (z-1)-z)}{p^2}, &\text{ otherwise. }
		\end{cases} 
	\end{align*}
	Hence the result follows.
\end{proof}

 As a corollary of the above theorem we get the following result.

\begin{cor}
	Let $G$ be a non-abelian group of order $p^n$ with  $|Z(G)| = p^{n-2}$. Then  CNL-spectrum, CNSL-spectrum, CNL-energy and CNSL-energy of CCC-graph of $G$ are given by\\
	$\cnlspec(\Gamma_G)= \left\{(0)^{p+1}, \left(\frac{1}{p^2}(p^{n - 1} - p^{n - 2})(p^{n - 1} - p^{n - 2} - 2p)\right)^{\frac{(p + 1)}{p}\left(p^{n - 1} - p^{n - 2} - p\right)}\right\}$,
	\begin{align*}
		&\cnqspec(\Gamma_G) = \bigg\{ \left(\frac{2}{p^2}(p^{n - 1} - p^{n - 2} - p)(p^{n - 1} - p^{n - 2} - 2p)\right)^{p + 1},\\ &\qquad\qquad\qquad\qquad\qquad\qquad\qquad\qquad\left(\frac{1}{p^2}(p^{n - 1} - p^{n - 2} - 2p)^2\right)^{\frac{(p + 1)}{p}(p^{n - 1} - p^{n - 2} - p)} \bigg\}
	\end{align*}
	and
	\[
	LE_{CN}(\Gamma_G) = LE^+_{CN}(\Gamma_G) = \begin{cases}		
		0, \qquad\qquad\qquad\qquad\qquad\qquad\qquad\text{ for } p=2  \,\&\, n=3\vspace{0.3 cm}\\
		\frac{2 (p+1) \left(p \left(p^{n-2}-2\right)-p^{n-2}\right) \left(p \left(p^{n-2}-1\right)-p^{n-2}\right)}{p^2}, \qquad\text{ otherwise. }
	\end{cases} 
	\]  
\end{cor}
\begin{proof}
	Here $\frac{G}{Z(G)}\cong \mathbb{Z}_p \times \mathbb{Z}_p$. Hence the result follows from Theorem \ref{Z_pZ_p theorem}.
\end{proof}

\subsection{Groups whose central quotient is isomorphic to a dihedral group}
The CCC-graph of this class of group was first studied by Salahshour \cite{Salah-2020} in 2020. Salahshour \cite[Theorem 1.2]{Salah-2020} obtained the following structures of $\Gamma_G$ (where $\frac{G}{Z(G)} \cong D_{2n}$)
\begin{equation}\label{eq-ccc-G/Z(G) D_2m}
\Gamma_G = \begin{cases}
	K_{\frac{(n-1)|Z(G)|}{2}}\cup 2K_{\frac{|Z(G)|}{2}}, &\text{ if } 2\mid n\\
	K_{\frac{(n-1)|Z(G)|}{2}}\cup K_{|Z(G)|}, &\text{ if } 2\nmid n.
\end{cases} 
\end{equation}
The CN-energy of such $\Gamma_G$ was studied in \cite{JN-2025}. Also spectrum, L-spectrum and Q-spectrum of commuting and non-commuting graphs of this class of groups was studied in \cite{DBN-2020,PJN-2018,DN-2017, DN-2018,  PDN-2018,  FN-2020} along with their respective energies.

\begin{thm}\label{P4_G_Z_D_2n theorem}
Let $G$ be a finite group such that $|Z(G)| = z$ and $\frac{G}{Z(G)}\cong D_{2n}$ (where $n \geq 3$). Then CNL-spectrum, CNSL-spectrum,  CNL-energy and CNSL-energy of CCC-graph of $G$ are as given below:
\begin{flushleft}
	${\rm(a)}$ If $n$ is even then\\ 
		{\rm (i)}\, $\cnlspec(\Gamma_G) = \left\lbrace(0)^3, (\frac{1}{4}(n z - z)(n z - z - 4))^{\frac{1}{2}(n z - z - 2)}, (\frac{1}{4} z (z - 4))^{z - 2}\right\rbrace$
		and \,
\[
LE_{CN}(\Gamma_G) = \frac{((n-1) z-2) (n (z+1) ((n-2) z-4)+11 z-4)}{2 (n+1)}.
\]

	{\rm (ii)} \, 
			$\cnqspec(\Gamma_G) = \left\lbrace (\frac{1}{2}(n z - z - 2)(n z - z - 4))^1, (\frac{1}{4} (n z - z - 4)^2)^{\frac{1}{2}(n z - z - 2)},\right.$\\			 			
			\qquad\qquad\qquad\qquad\qquad\qquad\qquad\qquad $\left. (\frac{1}{2}(z - 2)(z - 4))^2, (\frac{1}{4}(z - 4)^2)^{z - 2} \right\rbrace
		$\\
		and 
		\begin{align*}
			LE^+_{CN}(\Gamma_G) &= \begin{cases}
				\frac{28}{5}, &\text{ for } n=4 \,\&\, z=2\\
				\frac{3}{5} z^2 (4 z-6), &\text{ for } n=4 \,\&\, z \geq 3\\
				\frac{(n-2) (n-1) z^2 (n z-6)}{2 (n+1)}, &\text{ otherwise. }
			\end{cases}					 
		\end{align*}

	${\rm(b)}$ If $n$ is odd then\\
		{\rm (i)}\, $\cnlspec(\Gamma_G) = \left\lbrace (0)^2, (\frac{1}{4}(n z - z)(n z - z - 4))^{\frac{1}{2}(n z - z - 2)}, (z (z - 2))^{z - 1} \right\rbrace$\\
		and $LE_{CN}(\Gamma_G) = \begin{cases}
			0, &\text{ for } n=3 \,\&\, z=1\\
			 4(z-1) (z-2), &\text{ for } n=3 \,\&\, z\geq 2\\
			\frac{((n-1) z-2) \left((n-3) (n+1) z^2+((n-6) n+17) z-4 (n+1)\right)}{2 (n+1)}, &\text{ otherwise.}
		\end{cases}$\\
		{\rm (ii)}\, $\cnqspec(\Gamma_G) = \left\lbrace (\frac{1}{2}(n z - z - 2)(n z - z - 4))^1, (\frac{1}{4} (n z - z - 4)^2)^{\frac{1}{2}(n z - z - 2)},\right.$
\\		
		\qquad\qquad\qquad\qquad\qquad\qquad\qquad\qquad\qquad $\left.(2(z - 1)(z - 2))^1, ((z - 2)^2)^{z - 1} \right\rbrace$\\
		and $LE^+_{CN}(\Gamma_G) = \begin{cases}
			0, &\text{ for } n=3 \,\&\, z=1;~ n=5 \,\&\, z=1\\
			4(z - 1)(z - 2), &\text{ for } n=3 \,\&\, z\geq 2\\
			\frac{(n-5) (n-3) (n+3)}{2(n+1)}, &\text{ for } n\geq 7 \,\&\, z=1\\
			\frac{(n-3) (n-1) z^2 (n z+z-6)}{2 (n+1)}, &\text{ otherwise.}
		\end{cases}$
\end{flushleft}			
\end{thm}
\begin{proof}
From \eqref{eq-ccc-G/Z(G) D_2m}, we have $\Gamma_G = 	K_{\frac{(n-1)z}{2}}\cup 2K_{\frac{z}{2}}$ or $K_{\frac{(n-1)z}{2}}\cup K_{z}$ according as $n$ is even or odd.

(a){\rm (i)} If $n$ is even, then by Theorem \ref{CN-LE-CNSL-LE^+_Kn}
\small{
\[
\cnlspec(\Gamma_G) \!=\! \left\{(0)^1, \left(\frac{(n - 1) z}{2}\left(\frac{(n - 1)z}{2} - 2\right)\right)^{\frac{(n - 1)z}{2} - 1}, (0)^2, \left(\frac{z}{2}  \left(\frac{z}{2} - 2\right)\right)^{2(\frac{z}{2} - 1)}\right\}.
\]}	
\par Here $|V(\Gamma_G)| = \frac{1}{2} (n+1) z$ and $\tr(\cnrs(\Gamma_G)) = \frac{1}{8} z ((n ((n - 3) n + 3) + 1) z^2-6 ((n-2) n+3) z+8 (n+1))$. So, $\Delta(\Gamma_G) = \frac{(n ((n-3) n+3)+1) z^2-6 ((n-2) n+3) z+8 (n+1)}{4 (n+1)}$. Note that $z \geq 2$. We have
\begin{align*}
	L_1 &:= \left|0 - \Delta(\Gamma_G) \right| = \left| - \frac{(n ((n - 3) n + 3) + 1) z^2 - 6 ((n - 2) n + 3) z + 8 (n + 1)}{4 (n +  1)} \right|.
\end{align*}	
Let $\alpha_1(n, z) = (n ((n - 3) n + 3) + 1) z^2 - 6 ((n - 2) n + 3) z + 8 (n + 1)$. Then  $\alpha_1(n, z) = 8 + 8n + 6z(2n - 3) +  z^2 + 3n z^2 + \frac{n^2 z}{2}(nz - 12) + \frac{n^2 z^2}{2}(n - 6) > 0$    for $n \geq 12$. Also,  
$\alpha_1(4, z) = 29 z^2 - 66 z + 40 \geq 0$, $\alpha_1(6, z) = 127 z^2 - 162 z + 56 \geq 0$,  $\alpha_1(8, z) = 345 z^2 - 306 z + 72 \geq 0$ and $\alpha_1(10, z) = 731 z^2 - 498 z + 88 \geq 0$. Therefore
\[
L_1 = \frac{(n ((n - 3) n + 3) + 1) z^2 - 6 ((n - 2) n + 3) z + 8 (n + 1)}{4 (n +  1)}.
\]

We have
\begin{align*}
	L_2 &:= \left|\frac{1}{4} (n z - z) (n z - z - 4) - \Delta(\Gamma_G) \right| = \left| \frac{n (z + 1) ((n - 2) z - 4) + 11 z - 4}{2 (n + 1)} \right|.
\end{align*}
Let $\alpha_2(n, z) = \{ n (z + 1) ((n - 2) z - 4) + 11 z - 4 \}$. Then   $\alpha_2(n, z) > 0$ for $n \geq 6$, since $n - 2 \geq 4 \implies z(n - 2) - 4 \geq 0 \implies n(z + 1)(z(n - 2) - 4) \geq 0$. Also, 
$\alpha_2(4,z) = 8z^2 + 3z - 20 \geq 0$.  
Therefore	
\[
L_2 = \frac{n (z + 1) ((n - 2) z - 4) + 11 z - 4}{2 (n + 1)}.
\]

We have
\begin{align*}
L_3 := \left| \frac{1}{4} z (z - 4) - \Delta(\Gamma_G) \right| = \left| \frac{14z - 8 - 8n -  nz(2(z + 8) + n( - 6 + (n - 3)z))}{4 (n + 1)} \right|.
\end{align*}
Let $\alpha_3(n, z) = 14z - 8 - 8n -  nz(2(z + 8) + n( - 6 + (n - 3)z))$. For $n \geq 10$, $n((n - 3)z - 6) > 0$ and $2(z + 8) > 0$. So, $\alpha_3(n, z) < 0$ for all $n \geq 10$. Also, 
$\alpha_3(4,z) = - 24 z^2 + 46 z - 40 \leq 0$, $\alpha_3(6,z) = - 120 z^2 + 134 z - 56 \leq 0$ and $\alpha_3(8,z) = - 336 z^2 + 270 z - 72 \leq 0$. Therefore
\begin{align*}
L_3 &= - \frac{14z - 8 - 8n -  nz(2(z + 8) + n( - 6 + (n - 3)z))}{4 (n + 1)}.
\end{align*}
Hence, by \eqref{LEcn}, we get
\begin{align*}
	LE_{CN}(\Gamma_G) &= 3\times L_1 + \frac{1}{2}(nz - z - 2)\times L_2 + (z - 2)\times L_3\\
	& = \frac{((n-1) z-2) (n (z+1) ((n-2) z-4)+11 z-4)}{2 (n+1)}.
\end{align*}
(a){\rm (ii)} If $n$ is even, then by Theorem \ref{CN-LE-CNSL-LE^+_Kn}
\begin{align*}
	\cnqspec(\Gamma_G) = \left\lbrace \left( 2 \left( \frac{(n - 1)z}{2} - 1 \right)\left( \frac{(n - 1)z}{2} - 2 \right) \right)^1, \left(\left( \frac{(n - 1)z}{2} - 2 \right)^2 \right)^{ \frac{(n - 1)z}{2} - 1},\right.\\
	\left. \left(\frac{1}{2}(z - 2)(z - 4)\right)^2, \left(\frac{1}{4}(z - 4)^2\right)^{z - 2} \right\rbrace.
\end{align*}
	
Here $|V(\Gamma_G)| = \frac{1}{2} (n+1) z$ and $\tr(\cnrs(\Gamma_G) = \frac{1}{8} z ((n ((n-3) n+3)+1) z^2-6 ((n-2) n+3) z+8 (n+1))$. So, $\Delta(\Gamma_G) = \frac{(n ((n-3) n+3)+1) z^2-6 ((n-2) n+3) z+8 (n+1)}{4 (n+1)}$. Note that $z \geq 2$. We have
\begin{align*}
		B_1 &:= \left| \frac{1}{2} (n z - z - 2) (n z - z - 4) - \Delta(\Gamma_G) \right|\\
		&= \left| \frac{\left(n \left(n^2+n-5\right)+1\right) z^2-6 n (n+2) z+8 n+30 z+8}{4 (n+1)} \right|.
\end{align*} 
Let $\beta_1(n, z) = (n (n^2 + n - 5) + 1) z^2 - 6 n (n + 2) z + 8 n + 30 z + 8$. Then $\beta_1(n, z) = 8 + 8n + 30z + z^2 + n z^2(n - 5) + \frac{nz}{2}(n^2 z - 24) + \frac{n^2 z}{2}(nz - 12)$. Clearly for $n \geq 12$, $\beta_1(n,z) > 0$, as $n^2 z - 24 \geq 0$ and $nz - 12 \geq 0$. It can be seen that  
$\beta_1(4, z) = 61 z^2 - 114 z + 40 \geq 0$,  
$\beta_1(6, z) = 223 z^2 - 258 z + 56 \geq 0$, $\beta_1(8, z) = 537 z^2 - 450 z + 72 \geq 0 $ and $\beta_1(10, z) = 1051 z^2 - 690 z + 88 \geq 0$. Therefore
\[
B_1 = \frac{\left(n \left(n^2 + n - 5\right) + 1\right) z^2 - 6 n (n + 2) z + 8 n + 30 z + 8}{4 (n + 1)}.
\]
 
We have 
\begin{align*}
B_2 := \left| \frac{1}{4} (n z - z - 4)^2 - \Delta(\Gamma_G) \right| = \left| \frac{n ((n - 2) z^2 - (n + 6) z + 4) + 13 z + 4}{2 (n + 1)} \right|.
\end{align*}
Let $\beta_2(n, z) = n ((n - 2) z^2 - (n + 6) z + 4) + 13 z + 4$. Then $\beta_2(n, z)  = 4 + 4n + 13z + \frac{nz}{3}(nz - 18) + \frac{n^2 z}{3}(z - 3) + \frac{nz^2}{3}(n - 6) > 0$  for $n \geq 6$ and $z \geq 3$.  It can be seen that  
$\beta_2(4, z) = 8z^2 - 27z + 20 	=- 2$ or $\geq 0$ according as $z =2 $ or $z\geq 3$ 
and $\beta_2(n, 2) = 2n(n - 8) + 30 = - 2$ or $\geq 0$ according as $n =4 $ or $n\geq 6$. 
Therefore
\begin{align*}
B_2 = \begin{cases}
	\frac{1}{5} , &\text{ for } n=4 \,\&\, z=2\\
	\frac{n ((n - 2) z^2 - (n + 6) z + 4) + 13 z + 4}{2 (n + 1)}, &\text{ otherwise.}
\end{cases}
\end{align*}

We have
\begin{align*}
B_3 &:= \left| \frac{1}{2} (z - 2) (z - 4) - \Delta(\Gamma_G) \right|\\
	&= \left| \frac{ - \left(\left(n^3 - 3 n^2 + n - 1\right) z^2\right) + 6 (n - 4) n z + 8 n + 6 z +8 }{4 (n + 1)} \right|.
\end{align*}
Let $\beta_3(n, z)= - ((n^3 - 3 n^2 + n - 1) z^2) + 6 (n - 4) n z + 8 n + 6 z + 8$. Then $\beta_3(n, z) = 8(1 - nz) + 8n(1 - z) + 2z(3 - 4n) + z^2(1 - n) + \frac{n^2 z}{2}(12 - nz) + \frac{n^2 z^2}{2}(6 - n) < 0$  for $n \geq 12$. It can be seen that  
$\beta_3(4, z) = - 19z^2 + 6z + 40  \leq 0$,  $\beta_3(6, z) = - 113z^2 + 78z + 56 \leq 0$, 
$\beta_3(8, z) = - 327 z^2 + 198 z + 72 \leq 0 $ and $\beta_3(10, z) = - 709 z^2 + 366 z + 88 \leq 0 $. Therefore
\[
B_3 = - \frac{-\left(\left(n^3-3 n^2+n-1\right) z^2\right)+6 (n-4) n z+8 n+6 z+8}{4 (n+1)}.
\]

We have
\begin{align*}
B_4 &:= \left| \frac{1}{4} (z - 4)^2 - \Delta(\Gamma_G) \right|
	= \left| \frac{n (8 - z (n ((n - 3) z - 6) + 2 (z + 10))) + 10 z + 8}{4 (n + 1)} \right|.
\end{align*}
Let $\beta_4(n, z) = n (8 - z (n ((n - 3) z - 6) + 2 (z + 10))) + 10 z + 8$. Then  $\beta_4(n, z) = 8n - 10nz + 10z - 10nz + 8 -2nz^2 + \frac{n^2 z}{2}(12 - nz) + \frac{n^2 z^2}{2}(6 - n) < 0$  for $n \geq 12$. It can be seen that 

$\beta_4(4, z) = - 24 z^2 + 26 z + 40 \leq 0$, $\beta_4(6, z) = - 120 z^2 + 106 z + 56  \leq 0$,
$\beta_4(8, z) = - 336 z^2 + 234 z + 72 \leq 0$ and $\beta_4(10, z) = - 720 z^2 + 410 z + 88 \leq 0$. Therefore
\[
B_4 = - \frac{n (8 - z (n ((n - 3) z - 6) + 2 (z + 10))) + 10 z + 8}{4 (n + 1)}.
\]

Hence, by \eqref{LE+cn}, we get
\begin{align*}
LE^+_{CN}(\Gamma_G) &= 1\times B_1 + \frac{1}{2}(nz - z- 2)\times B_2 + 2\times B_3 + (z - 2)\times B_4\\
		&= \begin{cases}
		\frac{28}{5}, &\text{ for } n=4 \,\&\, z=2\\
		\frac{3}{5} z^2 (4 z-6), &\text{ for } n=4 \,\&\, z \geq 3\\
		\frac{(n-2) (n-1) z^2 (n z-6)}{2 (n+1)}, &\text{ otherwise. }
		\end{cases}					 
\end{align*}		
(b){\rm (i)} If $n$ is odd, then by Theorem \ref{CN-LE-CNSL-LE^+_Kn}	
\[
\cnlspec(\Gamma_G) = \left\{ (0)^1, \left(\frac{(n - 1)z}{2}\left(\frac{(n - 1)z}{2} - 2\right)\right)^{\frac{(n - 1)z}{2} - 1}, (0)^1, (z (z - 2))^{z - 1} \right\}.
\]
Here $|V(\Gamma_G)| = \frac{1}{2} (n+1) z$ and $\tr(\cnrs(\Gamma_G) = \frac{1}{8} z (n z+z-4) (((n-4) n+7) z-2 (n+1))$. So, $\Delta(\Gamma_G) = \frac{(n z+z-4) (((n-4) n+7) z-2 (n+1))}{4 (n+1)}$. We have
\begin{align*}
L'_1 := \left| 0 - \Delta(\Gamma_G) \right| = \left| - \frac{(n z + z - 4) (((n - 4) n + 7) z - 2 (n + 1))}{4 (n + 1)} \right|.
\end{align*}
Let $\alpha'_1(n, z) = (n z + z - 4) (((n - 4) n + 7) z - 2 (n +1 ))$. Then $\alpha'_1(n, z) = (nz + z - 4)(7z - 2 + \frac{nz}{2}(n - 8) + \frac{n}{2}(nz - 4)) > 0$ for  $n \geq 8$, since $z \geq 1$. Again 
$\alpha'_1(3, z) = 16 z^2 - 48 z + 32 \geq 0$, $\alpha'_1(5, z) = 72 z^2 - 120 z + 48 \geq 0$ and $\alpha'_1(7, z) = 224 z^2 - 240 z + 64 \geq 0$, as $z \geq 1$. Therefore
\[
L'_1 = \frac{(n z + z - 4) (((n - 4) n + 7) z - 2 (n + 1))}{4 (n + 1)}.
\]

We have	
\small{
\begin{align*}
L'_2 &:= \left| \frac{1}{4} (n z - z) (n z - z - 4) - \Delta(\Gamma_G) \right|
	= \left| \frac{n^2 z^2 + n^2 z - 2 n z^2 - 6 n z - 4 n - 3 z^2 + 17 z - 4}{2(n + 1)} \right|.
\end{align*}
}
Let $\alpha'_2(n,z) = n^2 z^2 + n^2 z - 2 n z^2 - 6 n z - 4 n - 3 z^2 + 17 z - 4$. Then $\alpha'_2(n,z) = 17z - 4 + \frac{n}{2}(nz - 8) + \frac{nz}{2}(n - 12) + \frac{z^2}{2}(n^2 - 6) + \frac{nz^2}{2}(n - 4) > 0$  for $n \geq 8$. It can be seen that    
$\alpha'_2(3, z) = 8 z - 16 =- 8$ or $\geq 0$ according as $z = 1$ or $z \geq 2$;  
$\alpha'_2(5, z) = 12 z^2 + 12 z - 24 \geq 0$ and $\alpha'_2(7, z) = 32 z^2 + 24 z - 32 \geq 0$, as $z \geq 1$. Therefore
\begin{align*}
L'_2 = \begin{cases}
	1, &\text{ for } n=3 \,\&\, z=1\\
	\frac{n^2 z^2 + n^2 z - 2 n z^2 - 6 n z - 4 n - 3 z^2 + 17 z - 4}{2(n + 1)}, &\text{ otherwise. }
	\end{cases}
\end{align*}

We have
\begin{align*}
L'_3 &:= \left| z (z - 2) - \Delta(\Gamma_G) \right| 
	= \left| \frac{ - n^3 z^2 + 3 n^2 z^2 + 6 n^2 z + n z^2 - 20 n z - 8 n - 3 z^2 + 22 z - 8}{4 (n + 1)} \right|.
\end{align*}
Let $\alpha'_3(n, z) = - n^3  z^2 + 3 n^2 z^2 + 6 n^2 z + n z^2 - 20 n z - 8 n - 3 z^2 + 22 z - 8$. Then $\alpha'_3(n, z) = - 8 - 8n - 2z(10n - 11) - 3z^2 - \frac{n^2 z}{3}(nz - 18) - \frac{nz^2}{3}(n^2 - 3) - \frac{n^2 z^2}{3}(n - 9) < 0$  for  $n \geq 19$. It can be seen that 
$\alpha'_3(3, z) = 16 z - 32 = - 16$ or $\geq 0$ according as $z = 1$ or $z\geq 2$, 
$\alpha'_3(5, z) = - 48 z^2 + 72 z - 48 \leq 0$, $\alpha'_3(7, z) = - 192 z^2 + 176 z - 64 \leq 0$, $\alpha'_3(9, z) = - 480 z^2 + 328 z - 80 \leq 0$, $\alpha'_3(11, z) = - 960 z^2 + 528 z - 96 \leq 0$, $\alpha'_3(13, z) = - 1680 z^2 + 776 z - 112 \leq 0 $, $\alpha'_3(15, z) = - 2688 z^2 + 1072 z - 128 \leq 0$ and $\alpha'_3(17, z) = - 4032 z^2 + 1416 z - 144 \leq 0$. Therefore
\begin{align*}
L'_3 = \begin{cases}
	z -2, &\text{ for } n=3 \,\&\, z\geq 2\\
	- \frac{ - n^3 z^2 + 3 n^2 z^2 + 6 n^2 z + n z^2 - 20 n z - 8 n - 3 z^2 + 22 z - 8}{4 (n + 1)}, &\text{ otherwise. }
	\end{cases}
\end{align*} 
Hence, by \eqref{LEcn}, we get
\begin{align*}
LE_{CN}(\Gamma_G) &= 2\times L'_1 + \frac{1}{2}(n z - z - 2)\times L'_2 + (z - 1)\times L'_3\\
	&= \begin{cases}
	0, &\text{ for } n=3 \,\&\, z=1\\
	4(z-1) (z-2), &\text{ for }n=3 \,\&\, z\geq 2\\
	\frac{((n-1) z-2) \left((n-3) (n+1) z^2+((n-6) n+17) z-4 (n+1)\right)}{2 (n+1)}, &\text{ otherwise.}
	\end{cases} 
\end{align*}
(b){\rm (ii)}	If $n$ is odd, then by Theorem \ref{CN-LE-CNSL-LE^+_Kn}
\begin{align*}
\cnqspec(\Gamma_G) =& \left\lbrace  \left( 2 \left( \frac{(n - 1)z}{2} - 1 \right)\left( \frac{(n - 1)z}{2} - 2 \right) \right)^1,\right.\\ &\left(\left( \frac{(n - 1)z}{2} - 2 \right)^2 \right)^{\frac{(n - 1)z}{2} - 1},
(2(z - 1)(z - 2))^1, \left((z - 2)^2\right)^{z - 1} \Bigg{\}}.
\end{align*}
\par Here $|V(\Gamma_G)| = \frac{1}{2} (n+1) z$ and $\tr(\cnrs(\Gamma_G) = \frac{1}{8} z (n z+z-4) (((n-4) n+7) z-2 (n+1))$. So, $\Delta(\Gamma_G) = \frac{(n z+z-4) (((n-4) n+7) z-2 (n+1))}{4 (n+1)}$. We have
\begin{align*}
B'_1 &:= \left| \frac{1}{2} (n z - z - 2) (n z - z - 4) - \Delta(\Gamma_G) \right|\\
	&= \left| \frac{n^3 z^2 + n^2 z^2 - 6 n^2 z - 5 n z^2 - 12 n z + 8 n - 5  z^2 + 42 z + 8}{4 (n + 1)} \right|.
\end{align*} 
Let $\beta'_1(n, z) = n^3 z^2 + n^2 z^2 - 6 n^2 z - 5 n z^2 - 12 n z + 8 n - 5  z^2 + 42 z + 8$. Then $\beta'_1(n, z) = 8 + 8n + 42z + z^2(n(n - 5) - 5) + nz(n(nz - 6) - 12)$. For $n \geq 9$ we have $nz - 6 \geq 3$ which gives $n(nz - 6) - 12 > 0$ and $n(n - 5) - 5 > 0$. Thus, $\beta'_1(n, z)> 0$. Again 
$\beta'_1(3, z) = 16 z^2 - 48 z + 32 \geq 0$, $\beta'_1(5, z) = 120 z^2 - 168 z + 48 \geq 0 $ and $\beta'_1(7, z) = 352 z^2 - 336 z + 64 \geq 0 $, as $z\geq 1$. Therefore
\begin{align*}
B'_1 =\frac{n^3 z^2 + n^2 z^2 - 6 n^2 z - 5 n z^2 - 12 n z + 8 n - 5  z^2 + 42 z + 8}{4 (n + 1)}.
\end{align*}

We have
\begin{align*}
B'_2 &:= \left| \frac{1}{4} (n z - z - 4)^2 - \Delta(\Gamma_G) \right|= \left| \frac{n^2 z^2 - n^2 z - 2 n z^2 - 6 n z + 4 n - 3 z^2 + 19 z + 4}{2(n + 1)} \right|.
\end{align*}
Let $\beta'_2(n,z) = n^2 z^2 - n^2 z - 2 n z^2 - 6 n z + 4 n - 3 z^2 + 19 z + 4$. Then $\beta'_2(n,z) = 4 + 4n + 19z + \frac{nz}{4}(nz - 24) + \frac{n^2 z}{4}(z - 4) + \frac{z^2}{4}(n^2 - 12) + \frac{nz^2}{4}(n - 8) > 0$  for $n \geq 9$ and $z\geq 5$. It can be seen that  
$\beta'_2(3,z) = 16 - 8z = 8$ or $\leq 0$ according as $z=1$ or $z\geq 2$;
$\beta'_2(5,z) = 12 z^2 - 36 z + 24 \geq 0$;
$\beta'_2(7,z) = 32 z^2 - 72 z + 32  = - 8$ or $\geq 0$ according as $z=1$ or $z\geq 2$;  
$\beta'_2(n,1) = 20 - 4n \geq 0$ or $< 0$ according as $n=1, 3, 5$ or  $n\geq 7$;
$\beta'_2(n,2) = 2n(n - 8) + 30 \geq 0$ for all $n$ and 
$\beta'_2(n,3) = 6 n^2 - 32 n + 34 = - 8$ or $\geq 0$ according as $n= 3$ or $n \ne 3$.
Therefore
\begin{align*}
B'_2 &= \begin{cases}
	z -2, &\text{for }n=3 \,\&\, z\geq 2\\
	\frac{n^2 z^2 - n^2 z - 2 n z^2 - 6 n z + 4 n - 3 z^2 + 19 z + 4}{2(n + 1)}, &\text{otherwise. }
	\end{cases}
\end{align*} 

We have  
\begin{align*}
B'_3 &\!:=\! \left| 2 (z - 1) (z - 2) - \Delta(\Gamma_G) \right|\!
	=\! \left| \frac{- n^3 z^2 + 3 n^2 z^2 + 6 n^2 z + 5 n z^2 - 36 n z + 8 n + z^2 + 6 z + 8}{4(n + 1)} \right|.
\end{align*}
Let $\beta'_3(n, z) = - n^3 z^2 + 3 n^2 z^2 + 6 n^2 z + 5 n z^2 - 36 n z + 8 n + z^2 + 6 z + 8$. Then $\beta'_3(n, z) = - 8nz + 8 - 8nz + 8n - 20nz + 6z - \frac{n^2 z}{4}(nz - 24) - \frac{z^2}{4}(n^3 - 4) - \frac{nz^2}{4}(n^2 - 5) - \frac{n^2 z^2}{4}(n - 12) < 0$  for $n \geq 25$. Again, 
$\beta'_3(3, z) = 16z(z - 3) + 32 \geq 0$, $\beta'_3(5, z) = - 24 z^2 - 24 z + 48 \leq 0$, $\beta'_3(7, z) = - 160 z^2 + 48 z + 64 \leq 0$, $\beta'_3(9, z) = - 440 z^2 + 168 z + 80 \leq 0$, $\beta'_3(11, z) = - 912 z^2 + 336 z + 96 \leq 0$, $\beta'_3(13, z) = - 1624 z^2 + 552 z + 112 \leq 0$, $\beta'_3(15, z) = - 2624 z^2 + 816 z + 128 \leq 0$, $\beta'_3(17, z) = - 3960 z^2 + 1128 z + 144 \leq 0$, $\beta'_3(19, z) = - 5680 z^2 + 1488 z + 160 \leq 0$, $\beta'_3(21, z) = - 7832 z^2 + 1896 z + 176 \leq 0$ and $\beta'_3(23, z) = - 10464 z^2 + 2352 z + 192 \leq 0$. Therefore,
\begin{align*}
B'_3 = \begin{cases}
	z(z-3) + 2, &\text{ for } n=3 \,\&\, z\geq 1\\
	- \frac{n^3 ( - z^2) + 3 n^2 z^2 + 6 n^2 z + 5 n z^2 - 36 n z + 8 n + z^2 + 6 z + 8}{4(n + 1)}, &\text{ for } n\geq 5 \,\&\, z\geq 1.
	\end{cases}
\end{align*}

We have
\begin{align*}
B'_4 &:= \left| (z - 2)^2 - \Delta(\Gamma_G) \right|= \left| \frac{ - n^3 z^2 + 3 n^2 z^2 + 6 n^2 z + n z^2 - 28 n z + 8 n - 3 z^2 + 14 z + 8}{4(n + 1)} \right|.
\end{align*}
Let $\beta'_4(n, z) = - n^3 z^2 + 3 n^2 z^2 + 6 n^2 z + n z^2 - 28 n z + 8 n - 3 z^2 + 14 z + 8$. Then $\beta'_4(n, z) = - 8n(z - 1) - 14z(n - 1) + (8 - 6nz - 3z^2) - \frac{n^2 z}{3}(nz - 18) - \frac{nz^2}{3}(n^2 - 3) - \frac{n^2 z^2}{3}(n - 9) < 0$  for $n \geq 19$. It can be seen that 
$\beta'_4(3, z) = 32 - 16z = 16$ or $\leq 0$ according as $z=1$ or $z\geq 2$; 
$\beta'_4(5, z) = 24z(1 - 2z) + 48 = 24$ or $\leq 0$ according as $z=1$ or  $z\geq 2$; 
$\beta'_4(7, z) = - 192 z^2 + 112 z + 64 \leq 0$; $\beta'_4(9,z) = - 480 z^2 + 248 z + 80 \leq 0$; $\beta'_4(11,z) = - 960 z^2 + 432 z + 96 \leq 0$; $\beta'_4(13,z) = - 1680 z^2 + 664 z + 112 \leq 0$; $\beta'_4(15,z) = - 2688 z^2 + 944 z + 128 \leq 0$ and $\beta'_4(17,z) = - 4032 z^2 + 1272 z + 144 \leq 0$, as $z\geq 1$. Therefore,
\begin{align*}
B'_4 = \begin{cases}
	\frac{ - n^3 z^2 + 3 n^2 z^2 + 6 n^2 z + n z^2 - 28 n z + 8 n - 3 z^2 + 14 z + 8}{4(n + 1)}, &\text{ for } n=3,5 \,\&\, z=1\\
	- \frac{- n^3 z^2 + 3 n^2 z^2 + 6 n^2 z + n z^2 - 28 n z + 8 n - 3 z^2 + 14 z + 8}{4(n + 1)}, &\text{ otherwise.}
	\end{cases}
\end{align*}
Hence, by \eqref{LE+cn}, we get
\begin{align*}
LE_{CN}^+(\Gamma_G) &= 1\times B'_1 + \frac{1}{2}(n z - z - 2)\times B'_2 + 1\times B'_3 + (z - 1)\times B'_4\\
&= \begin{cases}
	0, &\text{ for } n=3 \,\&\, z=1;~ n=5 \,\&\, z=1\\
	4(z - 1)(z - 2), &\text{ for } n=3 \,\&\, z\geq 2\\
	\frac{(n-5) (n-3) (n+3)}{2(n+1)}, &\text{ for } n\geq 7 \,\&\, z=1\\
	\frac{(n-3) (n-1) z^2 (n z+z-6)}{2 (n+1)}, &\text{ otherwise.}
	\end{cases} 
\end{align*}
Hence the result follows.
\end{proof}

 As a corollary of the above Theorem \ref{P4_G_Z_D_2n theorem}, we get the following results. 
\begin{cor}\label{Dihedral-theorem-P4}
	The CNL-spectrum, CNSL-spectrum, CNL-energy  and CNSL-energy of CCC-graph of the dihedral group $D_{2n}$ (where $n\geq 3$) are as given below:
	\begin{flushleft}
		${\rm(a)}$ If $n$ is odd then\\
		{\rm (i)} $\cnlspec(\Gamma_{D_{2n}}) = \{ (0)^2, (\frac{1}{4}(n - 1)(n - 5))^{\frac{1}{2}(n-3)}\}$ 		and 
		\[ 
		LE_{CN}(\Gamma_{D_{2n}})=  
			\frac{(n-5) (n-3) (n-1)}{n+1}.
	 		\]
		{\rm (ii)}\, $\cnqspec(\Gamma_{D_{2n}}) = \{ (0)^1, (\frac{1}{2}(n - 3)(n - 5))^1, (\frac{1}{4}(n - 5)^2)^{\frac{1}{2}(n - 3)} \}$ and 
		\[
		LE^+_{CN}(\Gamma_{D_{2n}}) =  \frac{(n - 5) (n - 3) (n + 3)}{2(n + 1)}.
			\]
		${\rm(b)}$ If $n$ is even then\\
		{\rm (i)}\, $\cnlspec(\Gamma_{D_{2n}}) = \{ (0)^3, (\frac{1}{4}(n - 2)(n - 6))^{\frac{1}{2}(n - 4)} \}$ and
		\[
		LE_{CN}(\Gamma_{D_{2n}}) = 
			\frac{3 (n - 6) (n - 4) (n - 2)}{2 (n + 2)}.
		\]
		{\rm (ii)}\, $\cnqspec(\Gamma_{D_{2n}}) = \{ (0)^2, (\frac{1}{2}(n - 4)(n - 6))^1, (\frac{1}{4}(n - 6)^2)^{\frac{1}{4}(n - 4)} \}$ and 
		\[
		LE^+_{CN}(\Gamma_{D_{2n}}) = \begin{cases}
			\frac{28}{5}, &\text{ for } n=8\\
			\frac{(n-6) (n-4) (n-2)}{n+2}, &\text{ for } n\ne 8.
		\end{cases}
		\]
	\end{flushleft}
\end{cor}
\begin{proof}
We know that 
$\frac{D_{2n}}{Z(D_{2n})}\cong 	D_{2\times \frac{n}{2}}$ or $D_{2n}$ according as $n$ is even or odd. Therefore, by Theorem \ref{P4_G_Z_D_2n theorem}, we get the required result.
\end{proof}

\begin{cor}\label{T_{4n}}
	The CNL-spectrum, CNSL-spectrum, CNL-energy and CNSL-energy of CCC-graph of the dicyclic group $T_{4n}$ (where $n\geq 2$) are as given below:
	\begin{enumerate}
		\item \, $\cnlspec(\Gamma_{T_{4n}}) = \{ (0)^3, ((n - 1)(n - 3))^{n - 2}\}$ and 
		\[ 
		LE_{CN}(\Gamma_{T_{4n}})= 
			\frac{6 (n-3) (n-2) (n-1)}{n+1}.
		\]
		\item \, $\cnqspec(\Gamma_{T_{4n}}) = \{ (0)^2, (2(n - 2)(n - 3))^1, ((n - 3)^2)^{n - 2}\}$ and 
		\[
		LE^+_{CN}(\Gamma_{T_{4n}}) = \begin{cases}
			\frac{28}{5}, &\text{ for }n = 4\\
			\frac{4 (n-3) (n-2) (n-1)}{n+1}, &\text{ for }n \ne 4.
		\end{cases}
		\]
	\end{enumerate}
\end{cor}
\begin{proof}
	We know that   
	$\frac{T_{4n}}{Z(T_{4n})}\cong D_{2n}$. Therefore, by Theorem \ref{P4_G_Z_D_2n theorem}, we get the required result. 
\end{proof}

\begin{cor}
	The CNL-spectrum, CNSL-spectrum, CNL-energy  and CNSL-energy of CCC-graph of the group $U_{6n} =\langle x,y:~x^{2n}=y^3=1,~x^{-1}yx=y^{-1}\rangle$ (where $n\geq 2$) are as given below:
	\begin{enumerate}
		\item $\cnlspec(\Gamma_{U_{6n}}) = \{ (0)^2, (n(n - 2))^{2(n - 1)} \} $
		and 
		$LE_{CN}(\Gamma_{U_{6n}}) = 4 (n - 2) (n - 1).$\\
		\item $\cnqspec(\Gamma_{U_{6n}}) = \{ (2(n - 1)(n - 2))^2, ((n - 2)^2)^{2(n - 1)} \}$
		and
		$$LE^+_{CN}(\Gamma_{U_{6n}}) = 4 (n - 2) (n - 1).$$
	\end{enumerate}
\end{cor}
\begin{proof}
	We know that 
	$\frac{U_{6n}}{Z(U_{6n})}=D_{2\times 3}$. Therefore, by Theorem \ref{P4_G_Z_D_2n theorem}, we get the required result. 
\end{proof}

\begin{cor}\label{U_{(m,n)}}
	The CNL-spectrum, CNSL-spectrum, CNL-energy  and CNSL-energy of CCC-graph of the group $U_{(n,m)}$ (where  $m\geq 2$ and $n\geq 2$) are as given below:
	\begin{flushleft}
		${\rm(a)}$ If $m$ is even then
		\begin{align*}
			{\rm (i)} &\cnlspec{(\Gamma_{U_{(n,m)}})} =\\
			&\quad \left\{ (0)^3, (n(n - 2))^{2(n - 1)},
			\left(\frac{1}{4}(nm - 2n)(nm - 2n - 4)\right)^{\frac{1}{2}(nm - 2n - 2)} \right\}
		\end{align*}	
		and 
		\begin{align*}
			LE_{CN}(\Gamma_{U_{(n,m)}})
			&= \begin{cases}
				4 (n - 1) (2 n - 1), &\text{ for }m=2 ~ \& ~ n\geq 2\\
				6(n - 2)(n - 1), &\text{ for }m=4 ~ \& ~ n\geq 2\\
				8 (n - 1) (2 n - 1), &\text{ for }m=6 ~ \& ~ n\geq 2\\
				\frac{((m-2) n-2) (m (2 n+1) ((m-4) n-4)+44 n-8)}{2 (m+2)}, &\text{ otherwise.}
			\end{cases}
		\end{align*}
		\begin{align*}
		 {\rm(ii)} &\cnqspec{(\Gamma_{U_{(n,m)}})} = \Bigg\{ (2(n - 1)(n - 2))^2, ((n - 2)^2)^{2(n - 1)}, \\  & \quad\left.
			\left(\frac{1}{2}(nm - 2n - 2)(nm - 2n - 4)\right)^1,
			\left(\frac{1}{4}(nm - 2n - 4)^2\right)^{\frac{1}{2}(nm - 2n - 2)} \right\rbrace
		\end{align*}
		and
		\begin{align*}
			LE^+_{CN}(\Gamma_{U_{(n,m)}}) 
			& = \begin{cases}
				4 (n-1) (2 n-1), &\text{ for }m=2 ~\& ~ n\geq 2\\
				6(n - 2)(n - 1), &\text{ for }m = 4 ~\& ~ n\geq 2\\
				8 n (2 n-3)+8, &\text{ for }m=6 ~\& ~ n\geq 2\\
				\frac{24}{5}n^2(4 n - 3), &\text{ for }m=8 ~\& ~n\geq 2\\
				\frac{(m-4) (m-2) n^2 (m n-6)}{m+2}, &\text{ otherwise.}
			\end{cases}
		\end{align*}
		${\rm(b)}$ If $m$ is odd then
		\begin{align*}
			{\rm(i)}\cnlspec&{(\Gamma_{U_{(n,m)}})} =\\
			&\quad \left\{ (0)^2, (n(n - 2))^{n - 1},
			\left(\frac{1}{4}(nm - n)(nm - n - 4)\right)^{\frac{1}{2}(nm - n - 2)}\right\}
		\end{align*}	
		and 
		\begin{align*}
			LE_{CN}(\Gamma_{U_{(n,m)}})
			&= \begin{cases}
				4(n - 1)(n - 2), \qquad\qquad\qquad\qquad\!\!\quad\qquad\text{ for }m=3 \,\&\, n\geq 2\\
				\frac{((m-1) n-2) \left((m-3) (m+1) n^2+((m-6) m+17) n-4 (m+1)\right)}{2 (m+1)}, \qquad\text{ otherwise. }
			\end{cases}
		\end{align*}
		\begin{align*}
			\qquad{\rm(ii)} \cnqspec&{(\Gamma_{U_{(n,m)}})} = \left\lbrace (2(n - 1)(n - 2))^1, ((n - 2)^2)^{n - 1},\right.\\  &\!\!\!\!\!\!\left.  \left(\frac{1}{2}(nm - n - 2)(nm - n - 4)\right)^1, \left(\frac{1}{4}(nm - n - 4)^2\right)^{\frac{1}{2}(nm - n - 2)}  \right\rbrace
		\end{align*}
		and
		\begin{align*}
			LE^+_{CN}(\Gamma_{U_{(n,m)}}) = \begin{cases}
				4(n - 1)(n - 2), &\text{ for }m=3\,\&\,n\geq 2\\
				\frac{(m-3) (m-1) n^2 (m n+n-6)}{2 (m+1)}, &\text{ otherwise. }
			\end{cases}
		\end{align*}
	\end{flushleft}
\end{cor}
\begin{proof}
		We know that 
	$\frac{U_{(n,m)}}{Z(U_{(n,m)})}$ is isomorphic to $D_{2\times \frac{m}{2}}$ or $D_{2m}$ according as $m$ is even or odd. Therefore, by Theorem \ref{P4_G_Z_D_2n theorem}, we get the required result. 
\end{proof}
\begin{cor}\label{SD_{8n}}
	The CNL-spectrum, CNSL-spectrum, CNL-energy  and CNSL-energy of CCC-graph of the group $SD_{8n}$ (where $n\geq 2$) are as given below:
	\begin{flushleft}
		${\rm(a)}$ If $n$ is even then \\
		{\rm(i)}\, $\cnlspec(\Gamma_{SD_{8n}}) = \left\{(0)^3, ((2n - 1)(2n - 3))^{2n - 2}\right\}$  and 
		\[
		LE_{CN}(\Gamma_{SD_{8n}}) = \frac{12 (n-1) (4 (n-2) n+3)}{2 n+1}.
		\]
		{\rm(ii)}\, $\cnqspec(\Gamma_{SD_{8n}}) = \left\{(0)^2, (2(2n - 2)(2n - 3))^1, ((2n - 3)^2)^{2n - 2}\right\} $ and
		\[
		LE^+_{CN}(\Gamma_{SD_{8n}}) =  \begin{cases}
			\frac{28}{5}, &\text{ for } n = 2\\
			\frac{8 (n-1) (2 n-3) (2 n-1)}{2 n+1},&\text{ for } n\geq 4.
		\end{cases}
		\]
		${\rm(b)}$ If $n$ is odd then \\
		{\rm(i)}\, $\cnlspec(\Gamma_{SD_{8n}}) = \left\{(0)^2, (8)^3, ((2n - 2)(2n - 4))^{2n - 3}\right\}$ and 
		\[ 
		LE_{CN}(\Gamma_{SD_{8n}}) = \begin{cases}
			24, &\text{ for }n = 3\\
			\frac{4 (2 n-3) (5 (n-3) n+4)}{n+1}, & \text{ for } n\geq 5.
		\end{cases}
		\]
		{\rm(ii)}\, $\cnqspec(\Gamma_{SD_{8n}}) = \left\{ (2(2n - 3)(2n - 4))^1, ((2n - 4)^2)^{2n - 3}, 12^1, 4^3 \right\}$ and 
		\[
		LE^+_{CN}(\Gamma_{SD_{8n}}) =\begin{cases}
			24, & \text{ for }n = 3\\
			\frac{16 (n-3) (n-1) (2 n-1)}{n+1},& \text{ for } n\geq 5. 
		\end{cases}
		\]
	\end{flushleft}				
\end{cor}
\begin{proof}
		We know that 
	$\frac{SD_{8n}}{Z(SD_{8n})}$ is isomorphic to $D_{2\times 2n}$ or $D_{2n}$ according  as $n$ is even or odd. Therefore, by Theorem \ref{P4_G_Z_D_2n theorem}, we get the required result. 
\end{proof}

\begin{cor}\label{V_{8n}}
	The CNL-spectrum, CNSL-spectrum, CNL-energy  and CNSL-energy of CCC-graph of the group $V_{8n}$ (where $n\geq 2$) are as given below:
	\begin{flushleft}
		${\rm(a)}$ If $n$ is even then\\
		{\rm(i)}\, $\cnlspec(\Gamma_{V_{8n}}) = \left\{ (0)^5, ((2n - 2)(2n - 4))^{2n - 3} \right\}$ and
		\[
		LE_{CN}(\Gamma_{V_{8n}}) = \frac{20 (n-2) (n-1) (2 n-3)}{n+1}.
		\]
		{\rm(ii)}\, $\cnqspec(\Gamma_{V_{8n}}) = \left\{ (0)^4, (2(2n - 3)(2n - 4))^1, ((2n - 4)^2)^{2n - 3} \right\}$ and
		\[
		LE^+_{CN}(\Gamma_{V_{8n}}) = \frac{16 (n-2) (n-1) (2 n-3)}{n+1}.
		\]
		${\rm(b)}$ If $n$ is odd then \\
		{\rm(i)}\, $\cnlspec(\Gamma_{V_{8n}}) =\left\{ (0)^3, ((2n - 1)(2n - 3))^{2n - 2} \right\}$ and
		\[
		LE_{CN}(\Gamma_{V_{8n}}) = \frac{12 (n-1) (4 (n-2) n+3)}{2 n+1}.
		\]
		{\rm(ii)}\, $\cnqspec(\Gamma_{V_{8n}}) = \left\{(0)^2, (2(2n - 2)(2n - 3))^1, ((2n - 3)^2)^{2n - 2}\right\} $ and 
		\[
		LE^+_{CN}(\Gamma_{V_{8n}}) = \frac{8 (n-1) (2 n-3) (2 n-1)}{2 n+1}.
		\]
	\end{flushleft}
\end{cor}
\begin{proof}
	(a) If $n$ is even then, by \cite[Proposition 2.4]{SA-CA-2020}, we have  $\Gamma_{V_{8n}} = K_{2n - 2} \cup 2 K_2$.
	
	{\rm(i)} By Theorem \ref{CN-LE-CNSL-LE^+_Kn}, we get				
	\[
	\cnlspec(\Gamma_{V_{8n}}) = \{ (0)^5, ((2n - 2)(2n - 4))^{2n - 3} \}.
	\]					
	Here $|V(\Gamma_{V_{8n}})| = 2 (n + 1)$ and $\tr(\cnrs(\Gamma_{V_{8n}})) = 4 (n-2) (n-1) (2 n-3)$. So, $\Delta(\Gamma_{V_{8n}}) = \frac{2 (n-2) (n-1) (2 n-3)}{n+1}$. We have
	\begin{align*}
		L_1 := \left| 0 - \Delta(\Gamma_{V_{8n}}) \right| = \left| - \frac{2 (n-2) (n-1) (2 n-3)}{n+1} \right| = \frac{2 (n-2) (n-1) (2 n-3)}{n+1},
	\end{align*}
	since $- 2 (n-2) (n-1) (2 n-3) < 0$, as $n \geq 2$, so $2n - 3 > 0$, $n - 2 \geq 0$ and $n - 1 > 0$. Also
	\begin{align*}
		L_2 := \left| (2 n-2) (2 n-4) - \Delta(\Gamma_{V_{8n}}) \right| &= \left| \frac{10 (n-2) (n-1)}{n+1} \right|\\ &= \frac{10 (n-2) (n-1)}{n+1}, \text{ as } n\geq 2.
	\end{align*}
	Therefore, by \eqref{LEcn}, we get
	\begin{align*}
		LE_{CN}(\Gamma_{V_{8n}}) = 5\times L_1 + (2n - 3)\times L_2 = \frac{20 (n-2) (n-1) (2 n-3)}{n+1}.
	\end{align*}
	{\rm(ii)} By Theorem \ref{CN-LE-CNSL-LE^+_Kn}, we get
	\[
	\cnqspec(\Gamma_{V_{8n}}) = \left\{ (0)^4, (2(2n - 3)(2n - 4))^1, ((2n - 4)^2)^{2n - 3} \right\}.
	\]
	
	We have
	\begin{align*}
		B_1 &:= \left| 0 - \Delta(\Gamma_{V_{8n}}) \right| = L_1,
	\end{align*}
	\begin{align*}
		B_2 := \left| 2(2n - 3)(2n - 4) - \Delta(\Gamma_{V_{8n}}) \right| &= \left| \frac{2 (n-2) (n+3) (2 n-3)}{n+1} \right|\\
		&= \frac{2 (n-2) (n+3) (2 n-3)}{n+1}, ~~~~\text{ as } n\geq 2,
	\end{align*}
	and
	\begin{align*}
		B_3 := \bigg| (2 n - 4)^2 & - \Delta(\Gamma_{V_{8n}}) \bigg| = \bigg| \frac{2 (n-2) (3 n-7)}{n+1} \bigg| = \frac{2 (n-2) (3 n-7)}{n+1}, 
	\end{align*}
	as $n\geq 2$, so $2 (n-2) (3 n-7)\geq 0$.
	Therefore, by \eqref{LE+cn}, we get
	\begin{align*}
		LE^+_{CN}(\Gamma_{V_{8n}}) &= 4\times B_1 + 1\times B_2 + (2n - 3)\times B_3 = \frac{16 (n - 2) (n - 1) (2 n - 3)}{n + 1}.
	\end{align*}

	(b) If $n$ is odd then, \cite[Proposition 2.4]{SA-CA-2020}, we have  $\Gamma_{V_{8n}} = K_{2n - 1} \cup 2K_1 = \Gamma_{D_{2 \times 4n}}$. Hence, the result follows from Corollary \ref{Dihedral-theorem-P4}.
\end{proof}

\section{Some consequences}
In this section, we discuss some consequences of the results obtained in Section 2. Looking at the CNL-spectrum and CNSL-spectrum of CCC-graphs of the groups considered in Section 2, we get the following result.
\begin{thm}
	Let $G$ be a finite non-abelian group with center $Z(G)$. Then the CCC-graph of $G$ is (CNSL) CNL-integral if
	\begin{enumerate}
		\item  $\frac{G}{Z(G)}\cong \mathbb{Z}_p\times \mathbb{Z}_p$.
		\item  $\frac{G}{Z(G)}\cong D_{2n}$.
		\item  $G$ is isomorphic to $D_{2n}$, $T_{4n}$, $U_{6n}$, $U_{(n, m)}$,  $SD_{8n}$ and  $V_{8n}$. 
	\end{enumerate}
\end{thm}

Now we shall determine whether CCC-graphs of these groups are (CNSL) CNL-hyperenergetic.
\begin{thm}
	Let $G$ be a finite non-abelian group  and $\frac{G}{Z(G)} \cong \mathbb{Z}_p \times \mathbb{Z}_p$. Then the CCC-graph of $G$ is not (CNSL) CNL-hyperenergetic.
\end{thm}
\begin{proof}
Let $|Z(G)| = z$. Then $z\geq 2$ and	 $|V(\Gamma_G)| = \frac{(p^2 - 1)z}{p}$. By Theorem \ref{Z_pZ_p theorem} and \eqref{LEcn-Kn} 
	\begin{align*}
		LE_{CN}(K_{|V(\Gamma_G)|})  - LE_{CN}(\Gamma_G)
		= \begin{cases}
			\frac{4 \left((p-2) p \left(2 p^2+p-2\right)+4\right)}{p^2}, &\text{ for } p \geq 2 \,\&\, z=2\\
			16, &\text{ for } p=2 \,\&\, z=3\\
			\frac{2 p^3 z^2-2 p^2 z^2-4 p^2-2 p z^2+2 z^2}{p}, &\text{ otherwise.}			
		\end{cases}		
	\end{align*} 
	Let $f_1(p) = 4 \left((p-2) p \left(2 p^2+p-2\right)+4\right)$ and $f_2(p, z) = 2 p^3 z^2-2 p^2 z^2-4 p^2-2 p z^2+2 z^2$, where $z \geq 3$. Then $f_1(p) > 0$. Also,
	 $f_2(p, z)= \frac{2}{3} (p-3) p^2 z^2+\frac{2}{3} p^2 \left(p z^2-6\right)+\frac{2}{3} \left(p^2-3\right) p z^2+2 z^2 > 0$ for $p\geq 3$.
	For $p=2$ we have $f_2(p, z) = 6z^2 - 16 >0$. Hence, $LE_{CN}(K_{|V(\Gamma_G)|})  - LE_{CN}(\Gamma_G) > 0$. 
By Theorem \ref{Z_pZ_p theorem}, we also have	$LE_{CN}(\Gamma_G) = LE^+_{CN}(\Gamma_G)$. Therefore,
$LE^+_{CN}(\Gamma_G) - LE^+_{CN}(K_{|V(\Gamma_G)|}) > 0$.
Hence the result follows.	
\end{proof}
An immediate corollary of the above theorem is given below. 
\begin{cor}
	Let $G$ be a non-abelian group of order $p^n$ and center $|Z(G)|=p^{n - 2}$. Then CCC-graph of $G$ is not (CNSL) CNL-hyperenergetic.
\end{cor}

\begin{thm}
Let $G$ be a finite non-abelian group and $\frac{G}{Z(G)} \cong D_{2n}$ (where $n \geq 3$). Then the CCC-graph of $G$ is 
	\begin{enumerate}
		\item CNL-borderenergetic if $n=3,11 ~ \& ~ z=1$.
		\item CNL-hyperenergetic except for  $n=4, 6 ~ \& ~ z=2$; $ n=4~ \& ~ z=3$; $n=4 ~ \& ~ z=4$; $n=3 ~ \& ~ z\geq 2$; $n=5,7,9 ~ \& ~ z=1$ and $n=5 ~ \& ~ z=2,3$.  
		\item CNSL-borderenergetic for $n=3 ~ \& ~ z=1$.
		\item CNSL-hyperenergetic except for $n=4 ~ \& ~ z=2,3,4,5$;  $n=6, 8 ~ \& ~ z=2$; $n=6 ~ \& ~ z=3$; $n=5 ~ \& ~ z=1$; $n=3 ~ \& ~ z\geq 2$; $n\geq 7$ ($n$ is odd) ~ $\&$ ~ $z=1$; $n=5, 7, 9 ~ \& ~ z=2$ and $n=5 ~ \& ~ z=3, 4$.
	\end{enumerate}
\end{thm}
\begin{proof}
	We have
	$$|V(\Gamma_G)| = \begin{cases}
		\frac{1}{2}(n + 1)z, &\text{ for $n$ is even}\\
		\frac{1}{2}(n + 1)z, &\text{ for $n$ is odd.} 
	\end{cases}$$
	
	\noindent \textbf{Case 1:}  $n$ is even 
	
	In this case $z \geq 2$. By Theorem \ref{SD_{8n}} and \eqref{LEcn-Kn}, we have
\begin{align*}
LE_{CN}&(K_{|V(\Gamma_G)|})  - LE_{CN}(\Gamma_G)\\
& = \frac{-n^3 z^3+3 n^2 z^3+12 n^2 z^2-2 n z^3-18 n z^2-24 n z+12 z^2+12 z}{2(n + 1)}.
\end{align*} 
	Let $f_1(n, z) = -n^3 z^3+3 n^2 z^3+12 n^2 z^2-2 n z^3-18 n z^2-24 n z+12 z^2+12 z$. Then 
	$f_1(n, z) = \frac{1}{2} n^2 z^3 (6-n) +\frac{1}{2} n^2 z^2 (24-n z)-2 n z^3 + 6 z^2 (2-3 n) +12 z (1-2 n) < 0$ for $n \geq 6$ and $z \geq 4$.  
	We have $f_1(n, 2) = 8n^2(9 - n) + 72 - 136n < 0$ for $n \geq 10$. 
	Also $f_1(4, 2) = 168$, $f_1(6, 2) = 120$ and $f_1(8, 2) = -504$. Therefore,
	$
	f_1(n, 2) > 0$ or $< 0$ according as $n=4, 6$ or $n \geq 8$.
	We have $f_1(n, 3) = 27n^2(7 - n) + 144 - 288n < 0$ for $n \geq 8$. 
 
	Also $f_1(4, 3) = 288$ and $f_1(6, 3) = - 612$. Therefore, 
	$
	f_1(4, 3)> 0$ or $<0$ according as $n=4$ or $n\geq 6$.

	Now we need to check for $n = 4$ and $z \geq 4$. We have $f_1(4, z) = 12z^2(11 - 2z) - 84z < 0$ for $z\geq 6$. 
	Also $f_1(4, 4) = 240$ and $f_1(4, 5) = -120$. Therefore,
	$
	f_1(4, z) > 0$ or $< 0$ according as $z=4$ or $z\geq 5$.
	Hence, $LE_{CN}(K_{|V(\Gamma_G)|})  - LE_{CN}(\Gamma_G) > 0$ for $n=4, 6 ~ \& ~ z=2$ and $n=4 ~ \& ~ z=3, 4$. Otherwise, $LE_{CN}(K_{|V(\Gamma_G)|})  - LE_{CN}(\Gamma_G) < 0$.

By Theorem \ref{SD_{8n}} and \eqref{LEcn-Kn}, we also have
	\begin{align*}
		LE^+_{CN}&(K_{|V(\Gamma_G)|})  - LE^+_{CN}(\Gamma_G)\\ 
		&= \begin{cases}
			\frac{92}{5},\qquad\qquad\qquad\qquad\qquad\qquad \text{ for }n = 4 \,\&\, z = 2\\
			\frac{1}{10} z \left(-24 z^2+161 z-150\right)+4,~\! \text{ for }n = 4 \,\&\, z \geq 3\\
			\frac{ - n^3 z^3 + n^3 z^2 + 3 n^2 z^3 + 9 n^2 z^2 - 6 n^2 z - 2 n z^3 - 15 n z^2 - 12 n z + 8 n + 13 z^2 - 6 z + 8 }{2(n + 1)},\quad \text{otherwise. }\\
		\end{cases}		
	\end{align*} 
	Let $f_2(z) = \frac{1}{10} z \left(-24 z^2+161 z-150\right)+4$ and  $f_3(n, z) = - n^3 z^3 + n^3 z^2 + 3 n^2 z^3 + 9 n^2 z^2 - 6 n^2 z - 2 n z^3 - 15 n z^2 - 12 n z + 8 n + 13 z^2 - 6 z + 8$.
		Then 
		$f_2(z)= \frac{1}{10}z \left(z\left( 161 - 24z\right) - 150 \right) + 4 > 0$ or $< 0$ according as $z = 3, 4, 5$ or $z \geq 6$. 
Also, $f_3(n, z)= \frac{1}{3} n^3 (3 - z) z^2 + \frac{1}{3} (9 - n) n^2 z^3 + \frac{1}{3} n^2 z^2 (27 - n z) - 6 n^2 z + z^2 (13 - 2 n z)+\left(8 n - 15 n z^2\right) + (8 - 12 n z) - 6 z < 0$ for $n \geq 10$ and $z \geq 3$. 

We have $f_3(n, 2) = 4n^2(12 - n) + 48 - 92n < 0$ for $n \geq 12$. Also, $f_3(6, 2) = 360$, $f_3(8, 2) = 336$ and $f_3(10, 2) = - 72$. Therefore,
$f_3(n, 2)  > 0$ or $< 0$ according as $n=6, 8$ or $n\geq 10$.

We have $f_3(6, z) = z^2(463 - 120z) - 294z + 56 < 0$ for $z\geq 4$. Also, $f_3(6, 3) = 101$. Therefore,
	$
	f_3(6, z) > 0$ or $< 0$ according as $z=3$ or $z\geq 4$.
	We have $f_3(8, z) = z^2(981 - 336z) - 486z + 72 < 0$ for $z\geq 3$. Therefore,
	$
	f_3(n, z)  > 0$ if  $n= 6, 8 \,\&\, z=2$ and $n=6 \,\&\, z=3$. Otherwise, $f_3(n, z)  < 0$.	Hence, $LE^+_{CN}(K_{|V(\Gamma_G)|})  - LE^+_{CN}(\Gamma_G)> 0$ if $n=4 ~ \& ~ z=2,3,4,5$; $n=6, 8 ~ \& ~ z=2$ and $n=6 ~ \& ~ z=3$. Otherwise,  $LE^+_{CN}(K_{|V(\Gamma_G)|})  - LE^+_{CN}(\Gamma_G) < 0$.

\vspace{.5cm}

\noindent\textbf{Case 2:}  $n$ is odd
	
By Theorem \ref{SD_{8n}} and \eqref{LEcn-Kn}, we have
\begin{align*}
	LE_{CN}&(K_{|V(\Gamma_G)|})  - LE_{CN}(\Gamma_G)\\ 
		&= \begin{cases}
			0, &\text{ for } n=3 \,\&\, z=1\\
			4z^2 - 4,&\text{ for } n=3 \,\&\, z\geq 2\\
			\frac{24 z - 24 n z + 12 z^2 - 24 n z^2 + 12 n^2 z^2 - 3 z^3 + n z^3 + 
				3 n^2 z^3 - n^3 z^3}{2(n + 1)}, &\text{ otherwise.}
		\end{cases}		
\end{align*}
	Clearly  $4z^2 - 4 > 0$ for $z\geq 2$. Let $f_4(n,z) = 24 z - 24 n z + 12 z^2 - 24 n z^2 + 12 n^2 z^2 - 3 z^3 + n z^3 + 
	3 n^2 z^3 - n^3 z^3$. Then  $f_4(n,z) = \frac{1}{3}n^2 z^3(9 - n) + \frac{1}{3}n z^3(3 - n^2) + \frac{1}{3}n^2 z^2(36 - nz) - 3z^3 + 12z^2(1 - 2n) + 24z(1 - n) < 0$  for $n\geq 9$ and $z\geq 4$. We have $f_4(n, 1) = n^2(15 - n) + 33 - 47n < 0$ for $n \geq 15$. Also, $f_4(5, 1) = 48$, $f_4(7, 1) = 96$, $f_4(9, 1) = 96$, $f_4(11, 1) = 0$ and $f_4(13, 1) = - 240$. Therefore
	\[ f_4(n, 1) \begin{cases}
		= 0, &\text{ for } n=11\\
		> 0, &\text{ for } n=5, 7, 9\\
		< 0, &\text{ for } n \geq 13.
	\end{cases}\]
	We have $f_4(n, 2) = 8n^2(9 - n) + 72 - 136n < 0$  for $n \geq 9$,  $f_4(5, 2) = 192$ and $f_4(7, 2) = -96$. Therefore,
	$f_4(n,2) > 0$ or $< 0$ according as $n=5$ or $n\geq 7$.
	
	We have $f_4(n,3) = 27n^2(7 - n) + 99 - 261n < 0$ for $n\geq 7$ and $f_4(5,3)= 144$. Therefore,
	$f_4(n,3) > 0$ or $< 0$ according as $n=5$ or $n\geq 7$.
	
	Again, we have $f_4(5,z) = 48z^2(4 - z) - 96z < 0$ and $f_4(7,z) = 48z^2(9 - 4z) - 144z < 0$ for $z\geq 4$. Therefore
	\[
	f_4(n, z) \begin{cases}
		= 0, &\text{ for } n=11 \,\&\, z=1\\
		> 0, &\text{ for } n=5, 7, 9 \,\&\, z=1; n=5 \,\&\, z=3\\
		< 0, &\text{ otherwise. }  
	\end{cases}
	\]
	Hence
	\begin{align*}
		LE_{CN}(K_{|V(\Gamma_G)|}) &- LE_{CN}(\Gamma_G)\\ 
		& \begin{cases}
			= 0, &\text{for } n=3, 11 \,\&\, z=1\\
			> 0, &\text{for } n=3 \,\&\, z\geq 2;\\
			& n=5,7,9 \,\&\, z=1;~ n=5 ~ \& ~ z=2,3\\
			< 0, &\text{otherwise.}
		\end{cases}		
	\end{align*}
	
By Theorem \ref{SD_{8n}} and \eqref{LEcn-Kn}, we also have
	\begin{align*}
		LE^+_{CN}&(K_{|V(\Gamma_G)|})  - LE^+_{CN}(\Gamma_G)\\ 
		&= \begin{cases}
			0, \qquad\qquad\qquad\qquad \text{ for } n=3 \,\&\, z=1\\
			4, \qquad\qquad\qquad\qquad \text{ for }n=5 \,\&\, z=1\\
			4(z^2 - 1), \qquad\qquad\quad\!\! \text{ for }n=3 \,\&\, z\geq 2\\
			\frac{n^2 + 4n - 21}{n + 1}, \qquad\qquad\quad\text{ for }n\geq 7 \,\&\, z=1\vspace{0.2 cm}\\
			\frac{- n^3 z^3 + n^3 z^2 + 3 n^2 z^3 + 9 n^2 z^2 - 6 n^2 z + n z^3 - 21 n z^2 - 12 n z + 8 n - 3 z^3 + 19 z^2 - 6 z + 8}{2(n + 1)}, \text{ otherwise.}
		\end{cases}		
	\end{align*}
	Clearly $4(z^2 - 1) > 0$ for $z \geq 2$ and $n^2 + 4n - 21 > 0$ for $n \geq 7$. Let $f_5(n, z) = - n^3 z^3 + n^3 z^2 + 3 n^2 z^3 + 9 n^2 z^2 - 6 n^2 z + n z^3 - 21 n z^2 - 12 n z + 8 n - 3 z^3 + 19 z^2 - 6 z + 8$. Then  $f_5(n, z)= \frac{1}{4} n^3 z^2(4 - z) + \frac{1}{4}n^2z^3(12 - n) + \frac{1}{4}n^2z^2(36 - nz) + \frac{1}{4}nz^3(4 - n^2) - 21nz^2 + (8n - 12nz) - 6n^2z + z^2(19 - 3z) + (8 - 6z) < 0$  for $n \geq 13$ and $z\geq 7$. We have $f_5(n, 2) = 4n^2(12 - n) + 48 - 92n < 0$  for $n\geq 13$. Again $f_5(5, 2) = 288$, $f_5(7, 2) = 384$, $f_5(9, 2) = 192$ and $f_5(11, 2) = -480$. Therefore,
	$
	f_5(n, 2) > 0$ or $< 0$ according as $n=5,7,9$ or  $n\geq 11$.
	
	We have $f_5(n, 3) = 18n^2(8 - n) + 80 - 190n < 0$ for $n\geq 9$. Again, $f_5(5, 3) = 480$ and $f_5(7, 3) = -368$. Therefore,
	$
	f_5(n, 3) > 0$ or $< 0$ according as $n=5$ or $n\geq 7$.
	
	We have $f_5(n, 4) = 24n^2(13 - 2n) + 96 - 312n < 0$ for $n\geq 7$. Also, $f_5(5, 4) = 336$. Therefore,
	$
	f_5(n, 4) > 0$ or $< 0$ according as $n=5$ or	$n\geq 7$.

	We have $f_5(n, 5) = 10n^2(57 - 10n) + 78 - 452n < 0$ for $n\geq 7$. Also, $f_5(5, 5) = - 432$. Therefore, $f_5(n, 5) < 0$.
	
	We have $f_5(n, 6) = 36n^2(26 - 5n) + 8 - 604n < 0$  for $n\geq 7$. Also, $f_5(5, 6) = -2112$. Therefore, $f_5(n, 6) < 0$.
	
	Now we shall  check for $n=5,7,9,11$ and $z\geq 7$. We have $f_5(5,z) = 24z^2(11 - 2z) + 48 - 216z < 0$, $f_5(7, z) = 16z^2(41 - 12z) + 64 - 384z < 0$, $f_5(9, z) = 8z^2(161 - 60z) + 80 - 600z < 0$ and $f_5(11, z) = 96z^2(23 - 10z) + 96 -864z < 0$, as $z \geq 7$. Thus, 
	$
	f_5(n, z) > 0$ if  $n=5, 7, 9$ \,\&\, $z=2$ and $n=5$ \,\&\, $z=3, 4$. Otherwise, $f_5(n, z) < 0$. 
	
	Hence
	\begin{align*}
		LE^+_{CN}&(K_{|V(\Gamma_G)|})  - LE^+_{CN}(\Gamma_G)\\ 
		&= \begin{cases}
			= 0, &\text{for } n=3 \,\&\, z=1\\
			> 0, &\text{for } n=5 ~ \& ~ z=1;~ n=3 ~ \& ~ z\geq 2;\\
			& n\geq 7 ~ \& ~ z=1;~ n=5, 7, 9 ~ \& ~ z=2;~ n=5 ~ \& ~ z=3, 4\\
			< 0, &\text{otherwise.}
		\end{cases}		
	\end{align*}
Hence the result follows.	
\end{proof}
\noindent As a corollary of the above theorem we get the following results.
\begin{cor}
	The CCC-graph of $U_{6n}$ $(n \geq 2)$ is not (CNSL) CNL-hyperenergetic.
\end{cor}

\begin{cor} \label{D_2n-hyperenergetic}
Let $G = D_{2n}, T_{4n}, U_{6n}$,  $SD_{8n}$ or $U_{(n, m)}$. Then
\begin{enumerate}
\item $\Gamma_G$ is CNL-borderenergetic if and only if $G = D_6$,  $D_{22}$ and $U_{(n, 2)}$ for $n \geq 2$.
\item $\Gamma_G$ is CNL-hyperenergetic if and only if $G =  D_{2n}$ for $n \geq 13$; $T_{4n}$ for $n \geq 7$;  $SD_{8n}$ for $n \geq 4$ and $U_{(n, m)}$ except for $m = 3 \,\&\, n\geq 2$, $m = 5 \,\&\, n = 2, 3$, $m = 4 \,\&\, n\geq 2$, $m = 8 \,\&\, n=2$ and $m = 6 \,\&\, n\geq 2$.
\item $\Gamma_G$ is CNSL-borderenergetic if and only if $G = D_6$ and $U_{(n, 2)}$ for $n \geq 2$.
\item $\Gamma_G$ is CNSL-hyperenergetic if and only if $G = D_{2n}$ for $n$ is even and $n \geq 20$; $T_{4n}$ for $n \geq 10$;  $SD_{8n}$ for $n \geq 6$ and $U_{(n, m)}$ except for $m = 3 \,\&\, n \geq 2$, $m = 5, 7, 9 \,\&\, n=2$, $m = 5 \,\&\, n = 3, 4$, $m = 4 \,\&\, n \geq 2$, $m = 8 \,\&\, n=2$, $m = 6 \,\&\, n\geq 2$ and $m = 10 \,\&\, n=2$.	
\end{enumerate}
\end{cor}

\begin{cor}
	The CCC-graph of $V_{8n}$ ($n \geq 2$) is CNL-hyperenergetic for $n \geq 6$ and CNSL-hyperenergetic for $n \geq 4$.	  
\end{cor}
\begin{proof}
\textbf{Case 1:} $n$ is even

We have $|V(\Gamma_{V_{8n}})| = (2n + 2)$. By Theorem \ref{V_{8n}} and \eqref{LEcn-Kn}, we get
\small{ 
\begin{align*}
	LE_{CN}&(K_{|V(\Gamma_{V_{8n}})|})  - LE_{CN}(\Gamma_{V_{8n}}) 
	&= \frac{120-32 (n-4) (n-2) n}{n+1}
	\begin{cases}
		> 0, &\text{ for } 2 \leq n \leq 4\\
		< 0, &\text{ for } n \geq 6.			 
	\end{cases}		
\end{align*}
}
\small{
\begin{align*}
	LE^+_{CN}&(K_{|V(\Gamma_{V_{8n}})|})  - LE^+_{CN}(\Gamma_{V_{8n}}) 
	&= -\frac{6 (n-1) (n (5 n-19)+16)}{n+1}
	\begin{cases}
		> 0, &\text{ for } n = 2\\
		< 0, &\text{ for } n \geq 4.			 
	\end{cases}		
\end{align*}
}
Thus, $\Gamma_{V_{8n}}$ is not CNL-hyperenergetic if $n = 2, 4$ and $\Gamma_{V_{8n}}$ is  CNL-hyperenergetic if $n \geq 6$. Also, it is   not CNSL-hyperenergetic if $n = 2$ and $\Gamma_{V_{8n}}$ is  CNSL-hyperenergetic if $n \geq 4$.

\textbf{Case 2:} $n$ is odd

 We have  $\Gamma_{V_{8n}} = K_{2n - 1} \cup 2K_1 = \Gamma_{D_{2 \times 4n}}$. Then, by Corollary \ref{D_2n-hyperenergetic}, we have that $\Gamma_{V_{8n}}$ is  CNL-hyperenergetic if $n \geq 4$ and CNSL-hyperenergetic if $n \geq 5$. Hence, the result follows.
\end{proof}

\subsection{Comparing various  CN-energies}
In this subsection, we compare various CN-energies of  CCC-graphs of the groups considered in Section 2. 
\begin{thm} \label{Comare-01}
	Let $G$ be a finite group such that $|Z(G)|=z \geq 2$ and $\frac{G}{Z(G)} \cong \mathbb{Z}_p\times \mathbb{Z}_p$. 
	If $p=2 \,\&\, z=3$ or $p\geq 3 \,\&\, z=2$ then 
	\[
	E_{CN}(\Gamma_G) < LE_{CN}(\Gamma_G) = LE^+_{CN}(\Gamma_G).
	\]
	For all other cases,
	$E_{CN}(\Gamma_G) = LE_{CN}(\Gamma_G) = LE^+_{CN}(\Gamma_G)$.
\end{thm}
\begin{proof} 
	In view of Theorem \ref{Z_pZ_p theorem}, it is sufficient to compare $E_{CN}(\Gamma_G)$ and $LE_{CN}(\Gamma_G)$.
	By   Theorem \ref{Z_pZ_p theorem} and \cite[Theorem 2.9]{JN-2025}, we have
	\[
	LE_{CN}(\Gamma_G) - E_{CN}(\Gamma_G) = \begin{cases}
		3, &\text{ for }p=2 \,\&\, z=3\\
		\frac{8 (p-2) (p+1)}{p^2}, &\text{ for }p\geq 2 \,\&\, z=2\\
		0, &\text{ otherwise.}
	\end{cases}
	\]
	Clearly, $8(p-2)(p+1) = 0$ or $> 0$ according as $p=2$ or $p > 2$.
	Hence, the result follows.
\end{proof}

As a corollary to Theorem \ref{Comare-01} we have the following result.
\begin{cor}
	Let $G$ be a  non-abelian $p$-group of order $p^n$ and $|Z(G)|=p^{n - 2}$, where $p$ is a prime and $n \geq 3$. Then $E_{CN}(\Gamma_G) = LE_{CN}(\Gamma_G) = LE^+_{CN}(\Gamma_G)$.
\end{cor}

\begin{thm}
Let $G$ be a finite group   and $\frac{G}{Z(G)} \cong D_{2n}$ $(n\geq 3)$. 
If $n=3\,\&\, z\geq 1$ or $n=5 \,\&\, z=1$ $($where $|Z(G)|=z)$ then 
\[
E_{CN}(\Gamma_G) = LE_{CN}(\Gamma_G) = LE^+_{CN}(\Gamma_G).
\]
For all other cases, $E_{CN}(\Gamma_G) < LE^+_{CN}(\Gamma_G) < LE_{CN}(\Gamma_G)$.
\end{thm}
\begin{proof}
\noindent \textbf{Case 1:}   $n$ is even
 
 In this case $z\geq 2$. By Theorem \ref{P4_G_Z_D_2n theorem} and \cite[Theorem 2.14]{JN-2025}, we have
\begin{align*}
LE^+_{CN}&(\Gamma_G) - E_{CN}(\Gamma_G) \\
		&= \begin{cases}
			\frac{8}{5}, &\text{ for }n=4 \,\&\, z=2\\
			\frac{z^2(24z-91)+150z-120}{10}, &\text{ for }n=4 \,\&\, z\geq 3\\
			\frac{(n-2) (n-1) n z^3-(n (n (n+5)-17)+15) z^2+6 (n+1)^2 z-24 (n+1)}{2(n+1)}, &\text{ otherwise.}
		\end{cases}
	\end{align*}
	Let $f_1(z) = z^2(24z-91)+150z-120$ and  $f_2(n, z) = (n-2) (n-1) n z^3-(n (n (n+5)-17)+15) z^2+6 (n+1)^2 z-24 (n+1)$. Then  $f_1(z) > 0$    for $z\geq 4$. Also, $f_1(3) = \frac{159}{10}$. Therefore, $f_1(z) > 0$ for $z\geq 3$. 
	It can be seen that   	
	$f_2(n, z) = \frac{1}{3} n^3 (z-3) z^2+\frac{1}{3} (n-9) n^2 z^3+\frac{1}{3} n^2 z^2 (n z-15)+6 \left(n^2 z-4\right)+2 n z^3+(17 n-15) z^2+12 n (z-2)+6 z > 0$ for $n\geq 10$ and $z\geq 3$. Also, $f_2(n,2) = 4(n-2)(n-3)^2 > 0$ for $n \geq 6$.
	We have $f_2(6, z) = z^2(120z - 309)+294z-168 > 0$  and $f_2(8, z) = z^2(336z - 711)+486z -216 > 0$ for $z\geq 3$. Therefore, $f_2(n, z)>0$ for $n \geq 6$ and $z \geq 2$. Hence
	$$LE^+_{CN}(\Gamma_G) - E_{CN}(\Gamma_G) > 0.$$

Again 
\begin{align*}
LE^+_{CN}(\Gamma_G) &- LE_{CN}(\Gamma_G) \\
	&= \begin{cases}
		-\frac{8}{5}, &\text{ for }n=4 \,\&\, z=2\\
		\frac{-z(29z-66) - 40}{10}, &\text{ for }n=4 \,\&\, z\geq 3\\
		-\frac{(n ((n-3) n+3)+1) z^2-6 ((n-2) n+3) z+8 (n+1)}{2 (n+1)}, &\text{ otherwise.}
	\end{cases}		
\end{align*}
Clearly, $-z(29z-66) - 40 < 0$ for $z\geq 3$. Let $f_3(n, z) = -((n ((n-3) n+3)+1) z^2-6 ((n-2) n+3) z+8 (n+1))$. Then $f_3(n, z) = -\frac{1}{2} (n-6) n^2 z^2-\frac{1}{2} n^2 z (n z-12)-3 n z^2-6 (2 n-3) z-8 n-z^2-8 < 0$ for  $n\geq 6$ and $z\geq 2$. Therefore
$$LE^+_{CN}(\Gamma_G) - LE_{CN}(\Gamma_G) < 0.$$ 
Hence, $E_{CN}(\Gamma_G) < LE^+_{CN}(\Gamma_G) < LE_{CN}(\Gamma_G)$.

\noindent	\textbf{Case 2:}  $n$ is odd

By Theorem \ref{P4_G_Z_D_2n theorem} and \cite[Theorem 2.14]{JN-2025}, we have
\begin{align*}
		&LE^+_{CN}(\Gamma_G) - E_{CN}(\Gamma_G) \\
		&= \begin{cases}
			0, &\text{ for }n=3, 5 \,\&\, z=1\\
			0, &\text{ for }n=3 \,\&\, z\geq 2\\
			\frac{(n - 5)(n - 3)}{n+1}, &\text{ for }n\geq 7 \,\&\, z=1\\
			\frac{z \left((n-3) (n-1) (n+1) z^2-(n (n (n+5)-21)+23) z+6 (n+1)^2\right)-16 (n+1)}{2 (n+1)}, &\text{ otherwise.}
		\end{cases}
	\end{align*}
	Clearly $(n - 5)(n - 3) > 0$ for $n\geq 7$. Let 
	$f_4(n, z) = z \left((n-3) (n-1) (n+1) z^2 -\right.$ $\left. (n (n (n+5)-21)+23) z+6 (n+1)^2\right)-16 (n+1)$, where $n \geq 5$ and $z \geq 2$.
	 Then $f_4(n, z) =\frac{1}{4} n^3 (z-4) z^2+\frac{1}{4} (n-12) n^2 z^3+\frac{1}{4} \left(n^2-4\right) n z^3+\frac{1}{4} n^2 z^2 (n z-20)+2 \left(3 n^2 z-8\right)+(21 n-23) z^2+4 n (3 z-4)+3 z^3+6 z > 0$ for $n \geq 13$ and $z\geq 4$. We have $f_4(n, 2) = 4(n - 2)(n - 3)^2 > 0$ as $n\geq 5$; $f_4(n, 3) = 18n^2(n - 6) + 182n - 124 > 0$ for $n\geq 7$; and  $f_4(5, 3) = 336$. Therefore, $f_4(n,3) > 0$.
	 Further, for $z\geq 4$ we have  $f_4(5, z) = 24z^2(2z - 7) + 216z - 96 > 0$;  $f_4(7, z)= 16z^2(12z - 29) + 384z - 128>0$; $f_4(9, z) = 8z^2(60z - 121) + 600z - 160 > 0$ and $f_4(11, z) = 192z^2(5z - 9) + 864z - 192 > 0$. Therefore, $f_4(n, z) > 0$. Thus
	$$LE^+_{CN}(\Gamma_G) - E_{CN}(\Gamma_G) \begin{cases}
		=0, &\text{ for }n=3\,\&\, z\geq 1;\, n=5\,\&\, z=1\\
		>0, &\text{ otherwise.}
	\end{cases}$$

Again
\begin{align*}
LE^+_{CN}(\Gamma_G) - LE_{CN}(\Gamma_G) 
	= \begin{cases}
		0, &\text{ for }n=3 \,\&\, z=1\\
		0, &\text{ for }n=3 \,\&\, z\geq 2\\
		0, &\text{ for }n=5 \,\&\, z=1\\
		-\frac{(n-5)^2 (n-3)}{2 (n+1)}, &\text{ for }n\geq 7 \,\&\, z=1\\
		-\frac{(n z+z-4) (((n-4) n+7) z-2 (n+1))}{2 (n+1)}, &\text{ otherwise. }
	\end{cases}
\end{align*}
Clearly, $-(n-5)^2 (n-3) < 0$ for $n\geq 7$. Let $f_5(n, z) = - (n z+z-4) (((n-4) n+7) z-2 (n+1))$. Then $f_5(n, z) = -\frac{1}{2} (n-6) n^2 z^2-\frac{1}{2} n^2 z (n z-12)-3 n z^2-6 (2 n-5) z-8 n-7 z^2-8 < 0$ for $n \geq 7$ and $z\geq 2$. For $z\geq 2$ we have  $f_5(5, z) = -24(z - 1)(3z - 2) < 0$. Therefore, $f_5(n, z) < 0$. Thus
$$LE^+_{CN}(\Gamma_G) - LE_{CN}(\Gamma_G) \begin{cases}
	=0, &\text{ for }n=3\,\&\, z\geq 1;\, n=5\,\&\, z=1\\
	<0, &\text{ otherwise.}
\end{cases}$$
Hence, $E_{CN}(\Gamma_G) = LE^+_{CN}(\Gamma_G) = LE_{CN}(\Gamma_G)$, if $n=3\,\&\, z\geq 1$ or $n=5\,\&\, z=1$. For all other cases, $E_{CN}(\Gamma_G) < LE^+_{CN}(\Gamma_G) < LE_{CN}(\Gamma_G)$. 
This completes the proof.
\end{proof}
We conclude this section with the following corollary.
\begin{cor}
Let $G = D_{2n}, T_{4n}, U_{6n}, SD_{8n}, V_{8n}$ or $U_{(n, m)}$. Then
\begin{enumerate}
\item $E_{CN}(\Gamma_G) = LE^+_{CN}(\Gamma_G) = LE_{CN}(\Gamma_G)$ if and only if $G = D_{6}, D_{8}, D_{10}, D_{12}, T_{8}$, $T_{12}$, $SD_{28}, V_{16}$, $U_{6n}$ for $n \geq 2$ and $U_{(n, m)}$ for $m = 3, 4, 6$ and $n \geq 2$. 
\item $E_{CN}(\Gamma_G) < LE^+_{CN}(\Gamma_G) < LE_{CN}(\Gamma_G)$ if and only $G$ is not among the groups listed in {\rm (a)}.
\end{enumerate}
\end{cor}
In the following figures,  closeness of various CN-energies of CCC-graphs of $D_{2n}, T_{4n}$, $U_{6n}, SD_{8n}, V_{8n}$ and $U_{(n, m)}$  are depicted. 


\begin{minipage}[t]{.5\linewidth}
	\begin{tikzpicture}
		\begin{axis}
			[
			xlabel={$n$ (odd) $\rightarrow$},
			ylabel={CN-energies of $\Gamma_{D_{2n}}$ $\rightarrow$},
			xmin=3, xmax=20 ,
			ymin=0, ymax=150,
			grid = both,
			minor tick num = 1,
			major grid style = {lightgray},
			minor grid style = {lightgray!25},
			width=.7\textwidth,
			height=.7\textwidth,
			legend style={legend pos=north west},
			]
			\addplot[domain=2:30,samples at={3,5,7,...,21},mark=*,green, samples=24, mark size=.8pt]{(x-3)*(x-5)/2};
			\tiny
			\addlegendentry{$E_{CN}$}
			\addplot[domain=2:30,samples at={3,5,7,...,21},mark=triangle*,blue,mark size=.8pt, samples=24]{((x-5)*(x-3)*(x-1))/(x + 1)};
			\tiny
			\addlegendentry{$LE_{CN}$}
			\addplot[domain=2:30,samples at={3,5,7,...,21},mark=square*, red, mark size=.8pt, samples=24]{((x-5)*(x-3)*(x+3))/(2*(x + 1))};
			\tiny
			\addlegendentry{$LE^+_{CN}$}
		\end{axis}
	\end{tikzpicture}
	\vspace{-.2 cm}
	\captionsetup{font=footnotesize}
	
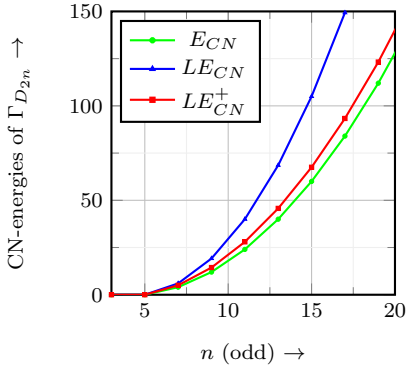
\captionof{figure}{CN-energies of $\Gamma_{D_{2n}}$, $n$ is odd}
\end{minipage}
\hspace{0.05cm}
\begin{minipage}[t]{.5\linewidth}
	\begin{tikzpicture}
		\begin{axis}
			[
			xlabel={$n$ (even) $\rightarrow$},
			ylabel={CN-energies of $\Gamma_{D_{2m}} \rightarrow$},
			xmin=2, xmax=20,
			ymin=0, ymax=150,
			grid = both,
			minor tick num = 1,
			major grid style = {lightgray},
			minor grid style = {lightgray!25},
			width=.7\textwidth,
			height=.7\textwidth,
			legend style={legend pos=north west},
			]
			\addplot[domain=10:50,samples at={4,6,8,10,12,14,16,18,20},mark=*,green, samples=20, mark size=.8pt]{(x-4)*(x-6)/2};
			\tiny
			\addlegendentry{$E_{CN}$}
			\addplot[domain=10:50,samples at={4,6,8,10,12,14,16,18,20},mark=triangle*,blue,mark size=.8pt, samples=20]{(3*(x-6)*(x-4)*(x-2))/(2*(x+2))};
			\tiny
			\addlegendentry{$LE_{CN}$}
			\addplot[domain=10:50,samples at={4,6,10,12,14,16,18,20},mark=square*, red, mark size=.8pt, samples=20]{((x-6)*(x-4)*(x-2))/(x+2)};
			\tiny
			\addlegendentry{$LE^+_{CN}$}
			\addplot[domain=7:9,samples at={8},mark=square*, red, mark size=.8pt, samples=20]{28/5};
			\tiny
		\end{axis}
	\end{tikzpicture}
	\vspace{-.2 cm}
	\captionsetup{font=footnotesize}
	
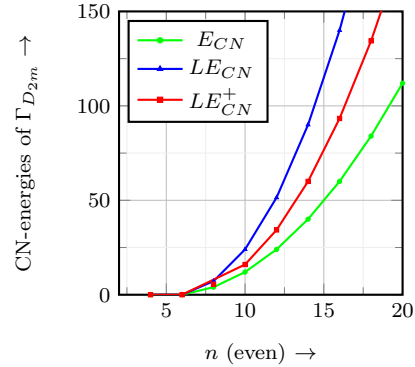
\captionof{figure}{CN-energies of $\Gamma_{D_{2n}}$, $n$ is even}
\end{minipage}

\vspace{.3cm}

\begin{minipage}[t]{.5\linewidth}
\begin{tikzpicture}
	\begin{axis}
		[
		xlabel={$n$ $\rightarrow$},
		ylabel={CN-ennergies of $\Gamma_{T_{4n}}$ $\rightarrow$},
		xmin=3, xmax=20,
		ymin=0, ymax=1500,
		grid = both,
		minor tick num = 1,
		major grid style = {lightgray},
		minor grid style = {lightgray!25},
		width=.7\textwidth,
		height=.7\textwidth,
		legend style={legend pos=north west},
		]
		\addplot[domain=2:21,samples at={2,3,4,5,6,7,8,9,10,11,12,13,14,15,16,17,18,19,20},mark=*,green, samples=24, mark size=.8pt]{2*(x-2)*(x-3)};
		\tiny
		\addlegendentry{$E_{CN}$}
		\addplot[domain=2:21,samples at={2,3,4,5,6,7,8,9,10,11,12,13,14,15,16,17,18,19,20},mark=triangle*,blue,mark size=.8pt, samples=24]{(6*(x-3)*(x-2)*(x-1))/(x+1)};
		\tiny
		\addlegendentry{$LE_{CN}$}
		\addplot[domain=2:21,samples at={2,3,5,6,7,8,9,10,11,12,13,14,15,16,17,18,19,20},mark=square*, red, mark size=.8pt, samples=24]{(4*(x-3)*(x-2)*(x-1))/(x+1)};
		\tiny
		\addlegendentry{$LE^+_{CN}$}
		\addplot[domain=3:5,samples at={4},mark=square*, red, mark size=.8pt, samples=20]{28/5};
	\end{axis}
\end{tikzpicture}
\vspace{-.2 cm}
\captionsetup{font=footnotesize}

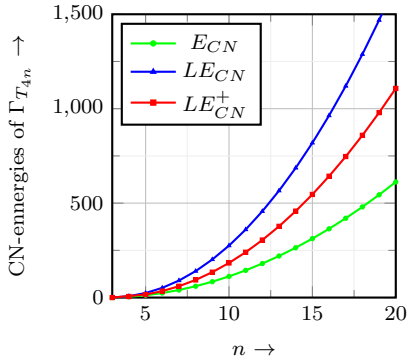
\captionof{figure}{CN-energies of $\Gamma_{T_{4n}}$}
\end{minipage}
\hspace{0.05cm}
\begin{minipage}[t]{.5\linewidth}
	\begin{tikzpicture}
		\begin{axis}
			[
			xlabel={$m$ (odd) $\rightarrow$},
			ylabel={CN-energies of $\Gamma_{U_{(4,m)}} \rightarrow$},
			xmin=2, xmax=20,
			ymin=0, ymax=7000,
			grid = both,
			minor tick num = 1,
			major grid style = {lightgray},
			minor grid style = {lightgray!25},
			width=.7\textwidth,
			height=.7\textwidth,
			legend style={legend pos=north west},
			]
			\addplot[domain=10:50,samples at={3,5,7,...,21},mark=*,green, samples=20, mark size=.8pt]{4*(2*x^2-7*x+9)};
			\tiny
			\addlegendentry{$E_{CN}$}
			\addplot[domain=10:50,samples at={5,7,9,...,21},mark=triangle*,blue,mark size=.8pt, samples=20]{(4*(2*x-3)*(5*x^2-15*x+4))/(x+1)};
			\tiny
			\addlegendentry{$LE_{CN}$}
			\addplot[domain=3:50,samples at={5,7,9,...,20},mark=square*, red, mark size=.8pt, samples=20]{(8*(x-3)*(x-1)*(4*x-2))/(x+1)};
			\tiny
			\addlegendentry{$LE^+_{CN}$}
			\addplot[domain=2:3,samples at={3},mark=square*, red, mark size=.8pt, samples=20]{24};
			\tiny
			\addplot[domain=2:50,samples at={3},mark=triangle*,blue,mark size=.8pt, samples=20]{24};
			\tiny
		\end{axis}
	\end{tikzpicture}
	\vspace{-.2 cm}
	\captionsetup{font=footnotesize}
	
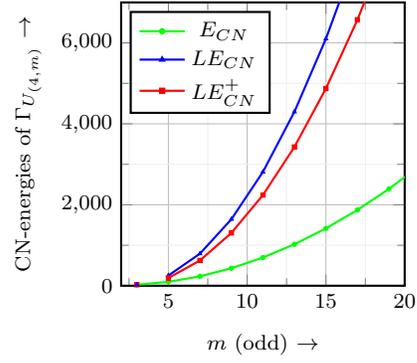
\captionof{figure}{CN-energies of $\Gamma_{U_{(4,m)}}$, $m$ is odd}
\end{minipage}
\vspace{0.3 cm}

\begin{minipage}[t]{.5\linewidth}
\begin{tikzpicture}
\begin{axis}
	[
	xlabel={$n$ (even) $\rightarrow$},
	ylabel={CN-energies of $\Gamma_{SD_{8n}}$ $\rightarrow$},
	xmin=3, xmax=20,
	ymin=0, ymax=5000,
	grid = both,
	minor tick num = 1,
	major grid style = {lightgray},
	minor grid style = {lightgray!25},
	width=.7\textwidth,
	height=.7\textwidth,
	legend style={legend pos=north west},
	]
	\addplot[domain=2:21,samples at={2,4,6,8,10,12,14,16,18,20},mark=*,green, samples=24, mark size=.8pt]{2*(2*x-2)*(2*x-3)};
	\tiny
	\addlegendentry{$E_{CN}$}
	\addplot[domain=2:21,samples at={2,4,6,8,10,12,14,16,18,20},mark=triangle*,blue,mark size=.8pt, samples=24]{(12*(x-1)*(4*(x-2)*x+3))/(2*x+1)};
	\tiny
	\addlegendentry{$LE_{CN}$}
	\addplot[domain=2:21,samples at={4,6,8,10,12,14,16,18,20},mark=square*, red, mark size=.8pt, samples=24]{(8*(x-1)*(2*x-3)*(2*x-1))/(2*x+1)};
	\tiny
	\addlegendentry{$LE^+_{CN}$}
	\addplot[domain=2:5,samples at={2},mark=square*, red, mark size=.8pt, samples=20]{28/5};
\end{axis}
\end{tikzpicture}
\vspace{-.2 cm}
\captionsetup{font=footnotesize}

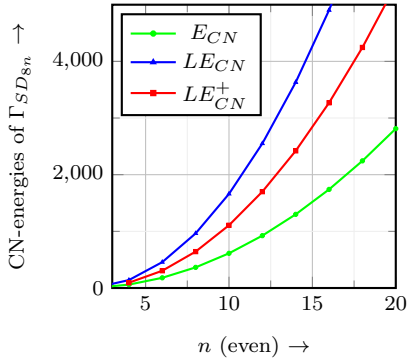
\captionof{figure}{CN-energies of $\Gamma_{SD_{8n}}$, $n$ is even}
\end{minipage}
\hspace{0.05cm}
\begin{minipage}[t]{.5\linewidth}
\begin{tikzpicture}
\begin{axis}
	[
	xlabel={$n$ (odd) $\rightarrow$},
	ylabel={CN-energies of $\Gamma_{SD_{8n}} \rightarrow$},
	xmin=3, xmax=20,
	ymin=0, ymax=6500,
	grid = both,
	minor tick num = 1,
	major grid style = {lightgray},
	minor grid style = {lightgray!25},
	width=.7\textwidth,
	height=.7\textwidth,
	legend style={legend pos=north west},
	]
	\addplot[domain=10:50,samples at={3,5,7,9,11,13,15,17,19,21},mark=*,green, samples=20, mark size=.8pt]{2*(2*x-3)*(2*x-4)+12};
	\tiny
	\addlegendentry{$E_{CN}$}
	\addplot[domain=10:50,samples at={5,7,9,11,13,15,17,19},mark=triangle*,blue,mark size=.8pt, samples=20]{(4*(2*x-3)*(5*(x-3)*x+4))/(x+1)};
	\tiny
	\addlegendentry{$LE_{CN}$}
	\addplot[domain=10:50,samples at={5,7,9,11,13,15,17,19},mark=square*, red, mark size=.8pt, samples=20]{(16*(x-3)*(x-1)*(2*x-1))/(x+1)};
	\tiny
	\addlegendentry{$LE^+_{CN}$}
	\addplot[domain=3:5,samples at={3},mark=square*, red, mark size=.8pt, samples=20]{24};
	\tiny
	\addplot[domain=3:5,samples at={3},mark=triangle*,blue, mark size=.8pt, samples=20]{24};
	\tiny
\end{axis}
\end{tikzpicture}
\vspace{-.2 cm}
\captionsetup{font=footnotesize}

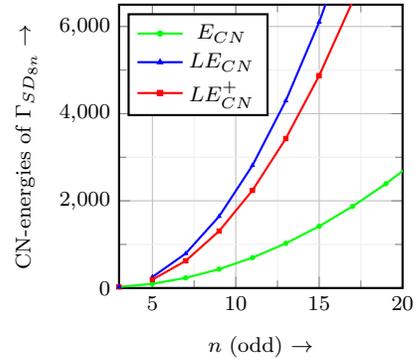
\captionof{figure}{CN-energies of $\Gamma_{SD_{8n}}$, $n$ is odd}
\end{minipage}

\begin{minipage}[t]{.5\linewidth}
\begin{tikzpicture}
\begin{axis}
[
xlabel={$n$ (even) $\rightarrow$},
ylabel={CN-energies of $\Gamma_{V_{8n}}$ $\rightarrow$},
xmin=3, xmax=20,
ymin=0, ymax=7000,
grid = both,
minor tick num = 1,
major grid style = {lightgray},
minor grid style = {lightgray!25},
width=.7\textwidth,
height=.7\textwidth,
legend style={legend pos=north west},
]
\addplot[domain=2:21,samples at={2,4,6,8,10,12,14,16,18,20},mark=*,green, samples=24, mark size=.8pt]{2*(2*x-3)*(2*x-4)};
\tiny
\addlegendentry{$E_{CN}$}
\addplot[domain=2:21,samples at={2,4,6,8,10,12,14,16,18,20},mark=triangle*,blue,mark size=.8pt, samples=24]{(20*(x-2)*(x-1)*(2*x-3))/(x+1)};
\tiny
\addlegendentry{$LE_{CN}$}
\addplot[domain=2:21,samples at={2,4,6,8,10,12,14,16,18,20},mark=square*, red, mark size=.8pt, samples=24]{(16*(x-2)*(x-1)*(2*x-3))/(x+1)};
\tiny
\addlegendentry{$LE^+_{CN}$}
\end{axis}
\end{tikzpicture}
\vspace{-.2 cm}
\captionsetup{font=footnotesize}

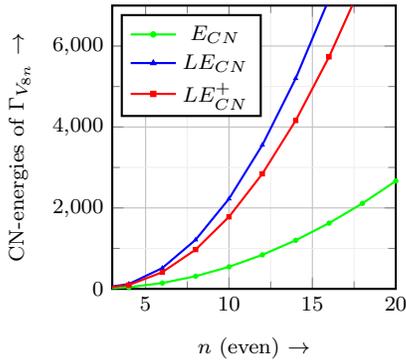
\captionof{figure}{CN-energies of $\Gamma_{V_{8n}}$, $n$ is even}
\end{minipage}
\hspace{0.05cm}
\begin{minipage}[t]{.5\linewidth}
\begin{tikzpicture}
\begin{axis}
[
xlabel={$n$ (odd) $\rightarrow$},
ylabel={CN-energies of $\Gamma_{V_{8n}} \rightarrow$},
xmin=2, xmax=20,
ymin=0, ymax=7000,
grid = both,
minor tick num = 1,
major grid style = {lightgray},
minor grid style = {lightgray!25},
width=.7\textwidth,
height=.7\textwidth,
legend style={legend pos=north west},
]
\addplot[domain=10:50,samples at={3,5,7,9,11,13,15,17,19,21},mark=*,green, samples=20, mark size=.8pt]{2*(2*x-2)*(2*x-3)};
\tiny
\addlegendentry{$E_{CN}$}
\addplot[domain=10:50,samples at={3,5,7,9,11,13,15,17,19,21},mark=triangle*,blue,mark size=.8pt, samples=20]{(12*(x-1)*(4*(x-2)*x+3))/(2*x+1)};
\tiny
\addlegendentry{$LE_{CN}$}
\addplot[domain=3:50,samples at={3,5,7,9,11,13,15,17,19,21},mark=square*, red, mark size=.8pt, samples=20]{(8*(x-1)*(2*x-3)*(2*x-1))/(2*x+1)};
\tiny
\addlegendentry{$LE^+_{CN}$}
\end{axis}
\end{tikzpicture}
\vspace{-.2 cm}
\captionsetup{font=footnotesize}

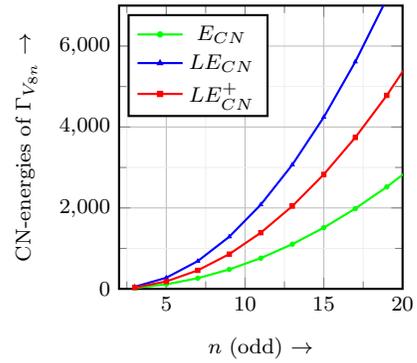
\captionof{figure}{CN-energies of $\Gamma_{V_{8n}}$, $n$ is odd}
\end{minipage}

\section{Conclusion}
In this paper, we compute common neighborhood (signless) Laplacian spectrum and energy of CCC-graphs of certain finite non-abelian groups.  We show that CCC-graphs of all the groups considered in this paper are CNL-integral and CNSL-integral.
The  common neighborhood spectrum and energy of CCC-graphs of theses groups are already computed in \cite{JN-2025}.   Analogous to the notion of super integral graph, we call a finite graph  \textit{super CN-integral} if it is CN-integral,  CNL-integral and CNSL-integral. Thus,  CCC-graphs of the groups considered in this paper are super CN-integral. It may be interesting to consider the following problem.
\begin{prob}
Characterize all finite non-abelian groups $G$ such that $\Gamma_G$ is super CN-integral. 
\end{prob}

The existence of finite non-abelian groups $G$ such that $\Gamma_G$ is CN-hyperenergetic is not clear (see \cite{JN-2025}). However, there are finite  non-abelian groups $G$ such that $\Gamma_G$ is CN-borderenergetic (See \cite[Theorem 3.6]{JN-2025}), CNL-hyperenergetic/CNL-borderenergetic and  CNSL-hyperenergetic/CSNL-borderenergetic (See Corollary \ref{D_2n-hyperenergetic}). Thus the following problem is worth considering.
\begin{prob}
Characterize all finite non-abelian groups $G$ such that 
$\Gamma_G$  is  CN-borderenergetic/ CNL-hyperenergetic/ CNL-borderenergetic/ CNSL-hyperenergetic/ \\CNL-borderenergetic. 
\end{prob}

We have found several classes of finite non-abelian groups $G$ such that $E_{CN}(\Gamma_G)  = LE_{CN}(\Gamma_G) = LE^+_{CN}(\Gamma_G)$ in Subsection 3.1. Thus, we pose the following problem.
\begin{prob}
	Characterize all finite non-abelian groups $G$ such that 
	$$E_{CN}(\Gamma_G)  = LE_{CN}(\Gamma_G) = LE^+_{CN}(\Gamma_G).$$ 
\end{prob}

In Subsection 3.1, we have also found  several classes of finite non-abelian groups $G$ such that $E_{CN}(\Gamma_G) < LE^+_{CN}(\Gamma_G) < LE_{CN}(\Gamma_G)$. In \cite[Theorem 4.6]{BN-2021}, it was observed that there are several classes of finite non-abelian groups $G$ such that  $E(\Gamma_G) < LE^+(\Gamma_G) < LE(\Gamma_G)$. It follows that there exist  finite non-abelian groups such that $E(\Gamma_G)$, $LE^+(\Gamma_G)$,  $LE(\Gamma_G)$  and $E_{CN}(\Gamma_G)$,  $LE^+_{CN}(\Gamma_G)$,  $LE_{CN}(\Gamma_G)$ behave similarly. Thus, the following problem arises naturally.
\begin{prob}\label{similar energies}
Determine all the finite non-abelian groups $G$ such that  $E(\Gamma_G)$, $LE^+(\Gamma_G)$,  $LE(\Gamma_G)$  and $E_{CN}(\Gamma_G),  LE^+_{CN}(\Gamma_G),  LE_{CN}(\Gamma_G)$ behave similarly. 
\end{prob}
It is worth noting that problem similar to Problem \ref{similar energies} can also be asked for any finite graph.	
	
In \cite{E-LE-Gutman}, Gutman et al. conjectured that $E(\mathcal{G}) \leq LE(\mathcal{G})$ for any finite graph $\mathcal{G}$ but soon after the announcement, this conjecture was refuted \cite{LL, SSM}. For the groups $G$, we consider in this paper, we have  
\begin{equation} \label{E-LE-cn}
E_{CN}(\Gamma_G) \leq LE_{CN}(\Gamma_G).
\end{equation}
In view of this it is too early to conjecture that the inequality \eqref{E-LE-cn} holds for CCC-graphs for any finite non-abelian group. However, one may consider the following problem.
\begin{prob}
Determine all the finite non-abelian groups such that the inequality \eqref{E-LE-cn} does not hold. In general, determine all the finite graphs $\mathcal{G}$ such that the inequality $E_{CN}(\mathcal{G}) \leq LE_{CN}(\mathcal{G})$ does not hold.
\end{prob}
\vskip 0.4 true cm

\begin{center}{\textbf{Acknowledgments}}
\end{center}
The first author would like to thank DST for the INSPIRE Fellowship (IF200226). 
\vskip 0.4 true cm

\end{document}